\documentclass[3p,12pt]{elsarticle}

\usepackage{amssymb}
\usepackage{amsmath}
\usepackage{natbib}
\setcitestyle{square,comma}
\usepackage{subcaption}
\usepackage{nomencl} 
\usepackage{romannum}
\usepackage{hyperref}
\usepackage{siunitx}
\makenomenclature

\usepackage{xcolor}

\usepackage{calligra,amsmath,amssymb}

\def\R{{\mathbb R}}

\def\wh{\widehat}

\AtBeginDocument{%
  \def\thefnote{\myfnsymbol{fnote}}}
\makeatletter
\def\myfnsymbol#1{\expandafter\@myfnsymbol\csname c@#1\endcsname}
\def\@myfnsymbol#1{%
  \ifcase #1
    \or $\dagger$%
    \or $\ddagger$%
    \else \@ctrerr
  \fi}
\def\fntext[#1]#2{\g@addto@macro\@fnotes{%
   \refstepcounter{fnote}\elsLabel{#1}%
   \def\thefootnote{\thefnote}%
   \global\setcounter{footnote}{\c@fnote}%
   \footnotetext{#2}}}
\makeatother
\usepackage{etoolbox}

\patchcmd{\dddot}{#1}{\kern0pt #1}{}{}

\makeatletter
\DeclareRobustCommand{\barredh}{%
  \mathord{\vphantom{h}\mathpalette\@barredh\relax}%
}

\newcommand{\@barredh}[2]{%
  \ooalign{%
    \hidewidth\@barredhbar#1\hidewidth\cr
    $\m@th#1h$\cr
  }%
}

\newcommand{\@barredhbar}[1]{%
  \check@mathfonts
  \ifx#1\displaystyle
    \fontsize{\f@size}{\z@}%
    \def\@barredhbarkern{0.3}%
  \else
    \ifx#1\textstyle
      \fontsize{\f@size}{\z@}
      \def\@barredhbarkern{0.3}%
    \else
      \ifx#1\scriptstyle
        \fontsize{\sf@size}{\z@}
        \def\@barredhbarkern{0.4}%
      \else
        \fontsize{\ssf@size}{\z@}
        \def\@barredhbarkern{0.47}%
      \fi
    \fi
  \fi
  \usefont{OT1}{cmr}{m}{n}%
  \kern-\@barredhbarkern em 
  \raisebox{-.175ex}{\symbol{'26}}%
}
\makeatother

\newcommand\bfu{{\mathbf u}}

\newcommand\bfI{{\mathbf I}}
\newcommand\bfA{{\mathbf A}}
\newcommand\bfB{{\mathbf B}}

\newcommand\bfF{{\mathbf F}}
\newcommand\bfG{{\mathbf G}}

\newcommand\bfK{{\mathbf K}}
\newcommand\bfL{{\mathbf L}}

\newcommand\bfP{{\mathbf P}}
\newcommand\bfQ{{\mathbf Q}}

\newcommand\bfS{{\mathbf S}}

\newcommand\bfU{{\mathbf U}}
\newcommand\bfV{{\mathbf V}}

\newcommand\bfX{{\mathbf X}}
\newcommand\bfY{{\mathbf Y}}
\newcommand\bfZ{{\mathbf Z}}

\newcommand\bfSigma{{\mathbf \Sigma}}

\def\phi{\varphi}

\journal{arXiv}

\begin{document}
\pagenumbering{arabic}
\begin{frontmatter}

\title{A high-order deterministic dynamical low-rank method for proton transport in heterogeneous media}

\author[1]{Pia Stammer\corref{cor1}}
\cortext[cor1]{p.k.stammer@tudelft.nl}
\author[2]{Niklas Wahl}
\author[3]{Jonas Kusch\fnref{fn1}}
\author[1]{Danny Lathouwers\fnref{fn1}}
\fntext[fn1]{Equal contributions, order alphabetical.}
\affiliation[1]{organization={Delft University of Technology, Dept. of Radiation Science and Technology}}
\affiliation[2]{organization={German Cancer Research Center -- DKFZ, Heidelberg Institute for Radiation Oncology (HIRO), National Center for Radiation Research in Oncology (NCRO)}}
\affiliation[3]{organization={Norwegian University of Life Sciences, Dept. of Data Science}}

\begin{abstract}
Dose calculations in proton therapy require the fast and accurate solution of a high-dimensional transport equation for a large number of (pencil) beams with different energies and directions. Deterministically solving this transport problem at a sufficient resolution can however be prohibitively expensive, especially due to highly forward peaked scattering of the protons. We propose using a model order reduction approach, the dynamical low-rank approximation (DLRA), which evolves the solution on the manifold of low-rank matrices in (pseudo-)time. For this, we compare a collided-uncollided split of the linear Boltzmann equation and its Fokker-Planck approximation. We treat the uncollided part using a ray-tracer and combine high-order phase space discretizations and a mixture model for materials with DLRA for the collided equation. 
Our method reproduces the results of a full-rank reference code at significantly lower rank, and thus computational cost and memory, and further makes computations feasible at much higher resolutions. At higher resolutions, we also achieve good accuracy with respect to TOPAS MC in homogeneous as well as heterogeneous materials. Finally, we demonstrate that several beam sources with different angles can be computed with little cost increase compared to individual beams.
\end{abstract}




\end{frontmatter}



\section{Introduction}\label{sec:1}
Accurate modeling of charged particle transport, particularly proton transport, is essential for applications such as radiation therapy, where dose deposition must be carefully controlled. The linear Boltzmann equation provides a physically accurate model for charged particle transport. However, solving it numerically is challenging due to multiscale effects, a high-dimensional phase space (space, angle, and energy), and strongly forward-peaked scattering. The latter necessitates a finely resolved phase space discretization or specialized methodological modifications, thus often leading to prohibitive memory requirements and computational costs. These computational requirements limit the feasibility of using the full Boltzmann model, particularly in settings that require many simulations, such as uncertainty quantification and treatment plan optimization in intensity-modulated proton therapy. 

For this reason, Monte Carlo methods, which can statistically solve the full physical model at high precision by simulating large numbers of particles, require acceleration techniques and strong computational resources to be practical \cite{vassiliev_feasibility_2008,elcim_dosimetric_2018}. Deterministic solutions to the Boltzmann equation remove the statistical noise and slow convergence of Monte Carlo methods and further allow access to full angle- and energy-resolved particle densities, thus better supporting further computations (e.g. uncertainty quantification) \cite{vassiliev_feasibility_2008,kawrakow_effect_2004}. Although the accuracy of such solvers can compete with that of Monte Carlo methods given equivalent physical modeling \cite{vassiliev_feasibility_2008,kan_review_2013}, they are not as wide-spread due to their high algorithmic complexity and challenges in efficient implementation. Rather, simplified models such as pencil beam algorithms, which apply dose kernels precomputed from Fermi-Eyges theory \cite{eyges_multiple_1948} or Gaussian fits to measurements and Monte Carlo simulations, are frequently used in practice \cite{hong_pencil_1996,schaffner_dose_1999,soukup_pencil_2005,yepes_comparison_2018}. While these models offer a computationally efficient alternative, they are based on simplifying assumptions which limit precision especially in heterogeneous media \cite{taylor_pencil_2017,yepes_comparison_2018}. Alternatively, deterministic approaches have been developed to efficiently resolve forward-peaked scattering, mainly based on multi-grid \cite{morel_angular_1991,turcksin_angular_2012} and adaptive grid refinement \cite{lathouwers_angular_2019,kophazi_spaceangle_2015}. While such refinement strategies can effectively resolve angular dependencies, they still require a sufficiently fine spatial resolution \cite{lathouwers_deterministic_2023}. Another approach in computational electron therapy \cite{kusch_robust_2023} reduces the memory footprint and computational costs by leveraging dynamical low-rank approximation (DLRA) \cite{koch_dynamical_2007}. This approach represents and evolves the solution to the linear Boltzmann equation in a low-rank format, thereby significantly reducing the number of degrees of freedom. Since DLRA is a dimensionality reduction method for time-dependent problems, a central idea in \cite{kusch_robust_2023} is to interpret the energy as a pseudo-time by using a continuous slowing down approximation \cite{larsen_electron_1997}. This interpretation treats energy loss in the medium as a surrogate for time evolution, enabling the use of time-dependent model reduction techniques such as DLRA.

A central challenge of DLRA is the inherent stiffness of the evolution equations of low-rank factors \cite{kieri_discretized_2016}. Therefore, robust time integrators that are not affected by this stiffness have been developed. The main idea of these integrators is to move only in flat subspaces of the low-rank manifold, thereby avoiding directions with a high curvature that result in stiff evolution equations \cite{lubich_projector-splitting_2014}. Perhaps the most popular class of robust time integrators is the (augmented) basis-update \& Galerkin (BUG) integrator \cite{ceruti_unconventional_2022,ceruti_rank-adaptive_2022}, which has since its introduction been modified to enable parallel update computations \cite{ceruti_parallel_2024} and to achieve high--order accuracy \cite{ceruti_robust_2024,kusch_second-order_2025}. Robust integrators for DLRA have been used extensively in kinetic theory~\cite{einkemmer_review_2025}, ranging from applications in plasma physics \cite{einkemmer_semi-lagrangian_2023,coughlin_robust_2024} to neutron criticality \cite{scalone_multi-fidelity_2025,kusch_low-rank_2022}. The use of DLRA in such problems relies heavily on the derivation of efficient spatial and temporal discretizations that preserve underlying structures.

While the method proposed in \cite{kusch_robust_2023} allows for highly efficient numerical simulations, it is currently limited to electron therapy, water-equivalent materials, as well as first-order temporal and spatial discretizations. An extension to proton simulations is challenging, primarily due to highly peaked scattering coefficients and the requirement of straggling terms in the continuous slowing down approximation. The latter makes the representation of energy as a pseudo-time more challenging and thus complicates the use of DLRA as a time-dependent model reduction technique. 

In this work, we propose a novel extension of DLRA to proton transport that addresses these challenges through the following contributions: (1) a ray-tracing model for uncollided particles to handle straggling, (2) high-order phase-space discretizations, and (3) a low-rank material decomposition. The raytracing model exploits the small number of relevant angular directions of uncollided particles to achieve low computational costs while being able to efficiently incorporate straggling effects. Collided particles are then advanced without straggling, thus enabling the use of dynamical low-rank approximation. Moreover, sharp gradients in proton transport require high-order phase-space discretizations. To mitigate computational costs imposed by limiters, we propose a nodal second--order upwind stencil combined with a fourth--order Runge--Kutta method in pseudo-time, i.e. energy, for numerical stability. We further use a modal spherical harmonics method (P$_{\text{N}}$) to discretize the directional domain. Lastly, to allow the incorporation of different  materials, we propose a low-rank material decomposition similar to \cite{kusch_low-rank_2022}, which enables efficient evaluation of the DLRA evolution equations.

The remaining manuscript is structured as follows: First we give a short introduction to proton transport and the dynamical low-rank approximation in Section~\ref{sec:background}. The former relies on the continuous slowing down approximation of the Boltzmann or Fokker-Planck equation, including a transformation that enables the use of energy as pseudo-time in the presence of material-dependent parameters. Then, the main method is introduced in Section~\ref{sec:methods}, including the ray tracer for the uncollided and spatial and angular discretizations as well as DLRA evolution equations based on the rank-adaptive BUG integrator for the time/energy update of the collided part of the equation. Finally, in Section~\ref{sec:numResults}, we present a variety of numerical experiments and in Section~\ref{sec:complexity}, we discuss the computational complexity and memory reduction. We provide a discussion and conclusion in Sections~\ref{sec:discussion} and \ref{sec:conclusion}.

\section{Background}\label{sec:background}
\subsection{Transport model for proton therapy} \label{sec:2}
We model the transport of particles through a medium using the continuous slowing down approximation of the steady-state linear Boltzmann equation
\begin{align}
    \mathbf{\Omega}\cdot\nabla_{\mathbf{r}}\psi(E,\mathbf{r},\mathbf{\Omega}) - \frac{\partial \mathcal{S}(E,\mathbf{r}) \psi(E,\mathbf{r},\mathbf{\Omega})}{\partial E} - \frac{1}{2}\frac{\partial^2 T(E,\mathbf{r})\psi(E,\mathbf{r},\mathbf{\Omega})}{\partial E^2}= \Gamma \psi(E,\mathbf{r},\mathbf{\Omega}),
\end{align}
where the phase space of the particle density $\psi$ consists of energy $E\in[E_{\mathrm{min}},E_{\mathrm{max}}]\subset \mathbb{R}_+$, space $\mathbf{r}\in \mathbb{R}^3$ and direction of flight $\mathbf{\Omega} \in \mathbb{S}^2$. We denote the stopping power by $\mathcal{S}$ and straggling by $T$. The collision operator $\Gamma$ on the right-hand side depends on the differential and total scattering cross sections $\Sigma_s$ and $\Sigma_t$ and is defined as
\begin{align}
    \Gamma \psi = \int_{\mathbb{S}^2} \Sigma_s(E,\mathbf{r},\mathbf{\Omega}'\cdot \mathbf{\Omega})\psi(E,\mathbf{r},\mathbf{\Omega}')\mathrm{d}\mathbf{\Omega}' - \Sigma_t(E,\mathbf{r})\psi(E,\mathbf{r},\mathbf{\Omega})
\end{align}
in case of the Boltzmann operator and as 
\begin{align}
    \Gamma\psi = \frac{\xi_1}{2}L
\end{align}
for the Fokker-Planck operator, where $\xi_1$ is the first transport coefficient defined as 
\begin{align}
 \xi_1 = 2\pi\int_{-1}^1 (1-\mu_0)\Sigma_s(E,\mathbf{r},\mu_0) \mathrm{d}\mu_0, \quad \mu_0 = cos(\theta) =  \mathbf{\Omega}'\cdot \mathbf{\Omega}
\end{align}
and $L$ is the spherical Laplace–Beltrami operator
\begin{align}
    L = \left[ \frac{\partial}{\partial \mu} (1-\mu^2) \frac{\partial}{\partial \mu}+\frac{1}{1-\mu^2}\frac{\partial^2}{\partial^2\Phi^2}\right],
\end{align}
where the spherical coordinates $\mu \in [-1,1]$ and $\Phi \in [0,2\pi) $ represent the cosine of the polar angle and the azimuthal angle, respectively.
Note, that, under the assumption of highly forward-peaked differential cross sections, the Fokker-Planck operator has been shown to be a first-order approximation to the Boltzmann operator \cite{leakeas_generalized_2001}. 

In order to later use energy as pseudo-time and separate angular and spatial dependency in the collision operator we apply the transformation $\widetilde{\psi}(E,\mathbf{r},\mathbf{\Omega}):= \mathcal{S}(E,\mathbf{r})\psi(E,\mathbf{r},\mathbf{\Omega})\;$
and assume that the tissue at each point in space is composed of a mixture of the same set  $\mathbb{A}=\{\text{H, C, N, O, Na, Mg, P, S, Cl, Ar,} \allowbreak \text{ K, Ca} \}$ of base materials with differing weights $w_i(\mathbf r), i \in \mathbb{A}$ \cite{burlacu_deterministic_2023}. Then the stopping power and cross sections can also be computed as weighted mixtures of the respective parameters for the individual elements using the Bragg additivity rule \cite{gottschalk_radiotherapy_2018}. Note that this model for material composition is much more accurate \cite{gottschalk_radiotherapy_2018, janni_proton_1982} than assuming water-equivalence as in \cite{kusch_robust_2023}. Then scattering cross sections $\Sigma_s$ are computed as 
\begin{align}
\label{eq:mixXS}
    \Sigma_s(E, \mathbf r, \mu_0) =  \sum_{i\in \mathbb{A}} w_i(\mathbf{r})\Sigma_{s;i}(E, \mu_0)\,,
\end{align}
where $\Sigma_{s;i}$ is the scattering cross-section of material $i$. The total scattering cross section $\Sigma_t(E, \mathbf r)$ is computed analogously, where for material $i$ we have $\Sigma_{t;i}(E) = 2\pi \int_{-1}^1 \Sigma_{s;i}(E, \mu_0) \mathrm d\mu_0 $. The resulting transport equation then takes the following form
\begin{align}
\label{eq:finalContinuous}
    -\frac{\partial  \widetilde{\psi}(E,\mathbf{r},\mathbf{\Omega})}{\partial E} &= - \mathbf{\Omega}\cdot\nabla_{\mathbf{r}} \frac{ \widetilde{\psi}(E,\mathbf{r},\mathbf{\Omega})}{\mathcal{S}(t,\mathbf{r})}  + \frac{1}{2}\frac{\partial^2 T(E,\mathbf{r}) \widetilde{\psi}(E,\mathbf{r},\mathbf{\Omega})}{\partial E^2 \mathcal{S}(E,\mathbf{r})} +  \sum_{i\in \mathbb{A}} w_i \Gamma_i  \frac{ \widetilde{\psi}(E,\mathbf{r},\mathbf{\Omega})}{\mathcal{S}(E,\mathbf{r})},\end{align}
    where for Boltzmann
  \begin{align}  \Gamma_{i} \psi &= \int_{\mathbb{S}^2} \Sigma_{s;i}(E,\mathbf{r},\mathbf{\Omega}'\cdot \mathbf{\Omega})\psi(E,\mathbf{r},\mathbf{\Omega}')\mathrm{d}\mathbf{\Omega}' - \Sigma_{t;i}(E,\mathbf{r})\psi(E,\mathbf{r},\mathbf{\Omega}) 
  \end{align}
  and for Fokker-Planck
  \begin{align}
     \Gamma_{i}\psi &= \frac{\xi_{1;i}}{2}L, \quad \xi_{1;i} = 2\pi\int_{-1}^1 (1-\mu_0)\Sigma_{s;i}(E,\mathbf{r},\mu_0)\mathrm{d}\mu_0 .
\end{align}

\subsubsection{Modeling of physical quantities}\label{sec:multiple_materials}
The particle and material dependent quantities in the transport equation include scattering cross sections, stopping power and straggling coefficients. The accuracy of our model strongly depends on a good modeling of these quantities. We determine the values for the set $\mathbb{A}$ containing 12 elements which are used to constitute tissue in the human body and mix these according to their proportions at each point in space. For comparability, we extract the stopping power for each material and different energies from the Monte Carlo code TOPAS MC \cite{perl_topas_2012}. Such data is not as easily obtained for straggling and differential cross sections, therefore we determine these based on physics models.

 We use William's model \cite{williams_passage_1997} to compute the modified straggling coefficient:

 \begin{equation}
     T(\mathbf r,E) = \sum_{i\in \mathbb{A}} \frac{1}{(4\pi\epsilon_0)^2} N_i(\mathbf r) 4 \pi e^4 Z_i\left( \frac{4I_i}{3m_ev_p(E)^2}\mathrm{ln}\frac{2m_ev_p(E)^2}{I_i} \right),
 \end{equation}
 where $\epsilon_0$ is the vacuum permittivity, $N_i(\mathbf r)$ are the atomic densities of atom $i$ within the material at position $\mathbf r$, $Z_i$ the respective atomic numbers, $I_i$ the ionization energies, $m_e$ the electron rest mass and $v_p(E)$ the velocity of the incident proton. 

For the differential elastic scattering, we use Molière's model \cite{scott_theory_1963,moliere_theorie_1948} which has been found in \cite{burlacu_yet_2024} to achieve good agreement with TOPAS MC and improved results compared to the small-angle first Born approximation. Here, the macroscopic differential cross section of the mixture material is given by
\begin{align}
\label{eq:elasticXS}
    \Sigma_s(E,\mathbf{r},\mu_0) =  \sum_{i\in \mathbb{A}} \tau_{i;\text{lab}}(\mu_0) \frac{4 N_i(\mathbf r) \alpha^2}{k(E)^2 (1-\mu_0 + \chi_{\alpha;i}(E))},
\end{align}
where $\chi_{\alpha;i}(E) = \chi_{0;i}(E)^2(1.13 + 3.76 a_i(E)^2)$ is the Molière corrected screening parameter according to \cite{scott_theory_1963,moliere_theorie_1948} with $\chi_{0;i}(E)=\frac{1.13\alpha Z_i^{1/3}m_ec}{p(E)}$ and $a_i = \frac{Z_i \alpha}{\beta(E)}$ representing the deviation from the Born approximation \cite{bethe_molieres_1953}. Further, $\beta=\frac{v_p(E)}{c}$ is the relative velocity, $\alpha$ the fine structure constant and $k(E)=\frac{p(E)}{\barredh}$ the reduced wave number with $p(E)$ the momentum and $\barredh$ the reduced Planck constant. Note, that equation \eqref{eq:elasticXS} contains a factor to convert the cross sections from the center of mass to the laboratory frame of reference, which can be computed as 
$$\tau_{i;\text{lab}}(\mu_0) = \frac{\left(1+2\mu_0\cdot\frac{m_p}{m_i}+\left(\frac{m_p}{m_i}\right)^2\right)^{\frac{3}{2}}}{1+\mu_0\cdot\frac{m_p}{m_i}},$$
where $m_p$ is the proton rest mass and $m_i$ the rest mass of the target nucleus \cite{burlacu_deterministic_2023}. Note, that since protons and neutrons have almost equal mass $\frac{m_p}{m_i}= A_i$ can be used to simplify the computation, where $A_i$ is the atomic mass number.

For simplicity, we neglect nuclear scattering (including production of secondary neutrons and recoil nuclei) and absorption as their effect on the deposited energy is relatively small \cite{newhauser_physics_2015,durante_nuclear_2016} and they do not affect the numerical method presented here.
 
\subsection{Dynamical low-rank approximation}
\label{sec:DLRA_intro}
Discretizing the transport equations presented in the previous section will result in a large matrix ordinary differential equation 
\begin{align}\label{eq:fullrank}
    \dot \bfu(t) = \bfF(t, \bfu(t))\, , \quad \bfu(t_0) = \bfu_0\,,
\end{align}
where $\bfu(t)\in\R^{n\times m}$ and energy will act as a pseudo-time. Here, $n$ refers to the number of spatial cells and $m$ is the number of moments of the modal angular discretization. 

A rank $r$ approximation to $\bfu(t)$ is given by $\bfu(t)\approx \bfX(t)\bfS(t)\bfV(t)^{\top}$ where $\bfX(t)\in\R^{n \times r}$ and $\bfV(t)\in\R^{m \times r}$ are orthonormal basis matrices and $\bfS(t)\in\mathbb{R}^{r\times r}$ is a not necessarily diagonal coefficient matrix. Here, the rank $r$ matrices form a smooth manifold, which we denote by $\mathcal{M}$. The core idea of dynamical low-rank approximation is to evolve the low-rank factors instead of the full solution in time \cite{koch_dynamical_2007}. This is achieved by projecting the dynamics of \eqref{eq:fullrank} onto the tangent space of $\mathcal{M}$ at the current solution. The projection restricts the dynamics of the solution to the low-rank manifold $\mathcal{M}$, thus ensuring a memory-efficient solution representation during the entire simulation. The projected equations then read
\begin{align}\label{eq:projectedDyn}
    \dot \bfu_r(t) = \bfP(\bfu_r(t))\bfF(t, \bfu_r(t))\,,
\end{align}
where for $\bfZ=\bfX\bfSigma\bfV^{\top}$, the projector onto the tangent space of $\mathcal{M}$ at $\bfZ$ takes the form 
\begin{align*}
    \bfP(\bfZ)\bfG = \bfX\bfX^{\top}\bfG(\bfI-\bfV\bfV^{\top}) + \bfG\bfV\bfV^{\top}\,.
\end{align*}

\begin{figure}[h!]
    \centering
    \includegraphics[width=0.5\linewidth]{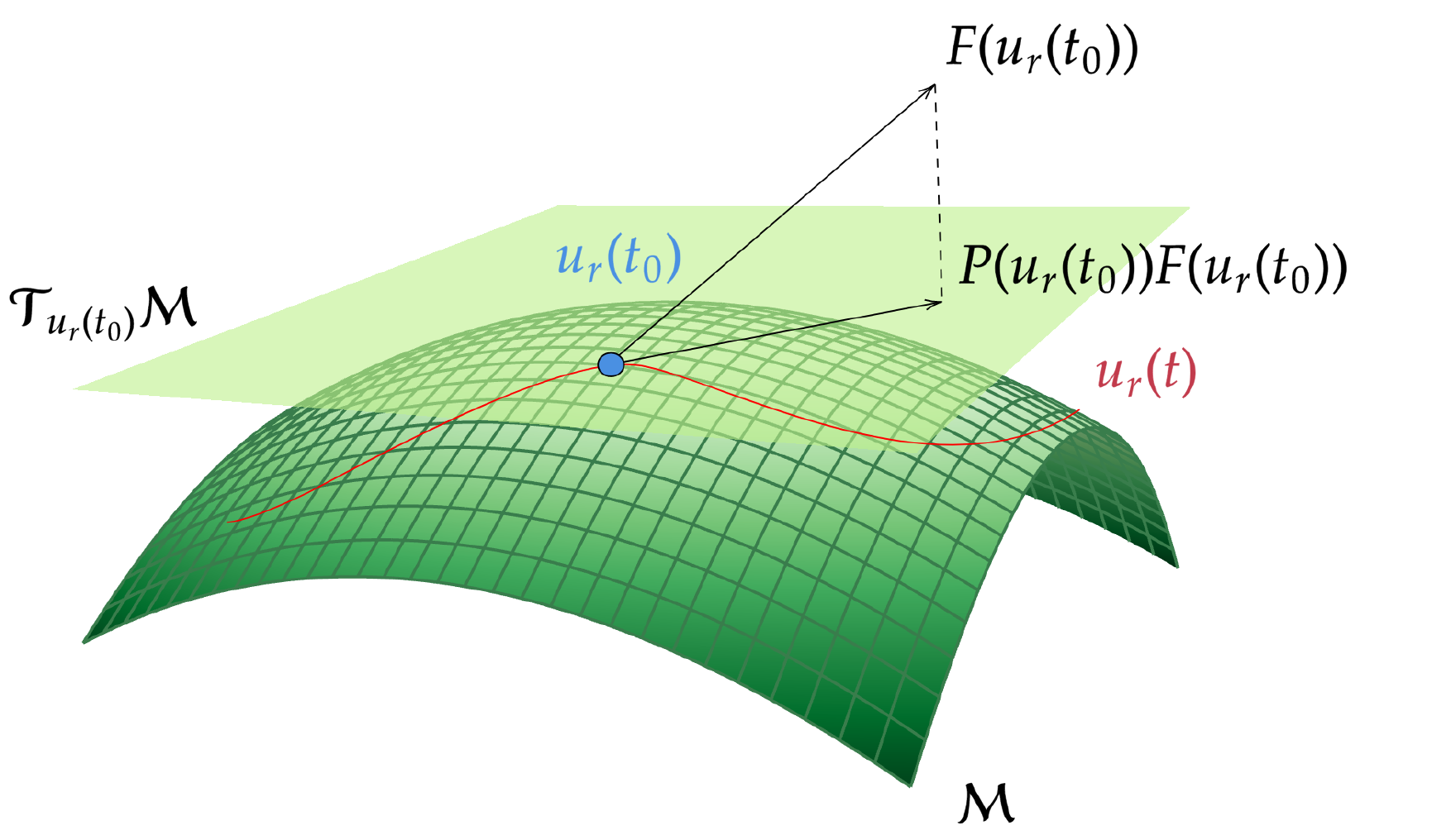}
    \caption{Graphical illustration of dynamical low-rank approximation. The solution at time $t_0$ lies on the manifold of rank $r$ matrices $\mathcal{M}$. It will not stay on $\mathcal{M}$ if it moves according to the original dynamics (i.e., into direction $\bfF(\bfu_r(t_0))$), but remains of rank $r$ if the original dynamics are projected onto the tangent space $\mathcal{T}_{\bfu_r(t)}\mathcal{M}$ via the projector $\bfP(\bfu_r(t))$.}
    \label{fig:manifold}
\end{figure}

A graphical illustration of dynamical low-rank approximation can be found in Figure~\ref{fig:manifold}. Since this manifold can exhibit high curvature, conventional time integration methods require prohibitively small step sizes. Therefore, careful construction of new time integration methods that take the geometry of $\mathcal{M}$ into account when solving \eqref{eq:projectedDyn} is required \cite{kieri_discretized_2016}. Robust time integration methods include the projector--splitting integrator~\cite{lubich_projector-splitting_2014} and basis-update \& Galerkin (BUG) integrators~\cite{ceruti_unconventional_2022,ceruti_rank-adaptive_2022,ceruti_parallel_2024}. The main idea of these integrators is to move on flat subspaces of the low-rank manifold, thereby avoiding directions with high curvature. The augmented BUG integrator \cite{ceruti_rank-adaptive_2022}, which we use in this work, first updates the basis matrices $\bf{U}$ and $\bfV$ in parallel, constraining the dynamics to
\begin{align}
    \{\bfK \bfV(t_0)^{\top} \,|\, \bfK \in\R^{n \times r}\} \subset \mathcal{M} \enskip\text{ and }\enskip \{\bfU(t_0)\bfL^{\top} \,|\, \bfL \in\R^{m \times r}\}\,\subset \mathcal{M}.
\end{align}
This leads to the differential equations
\begin{subequations}\label{eq:KLS}
\begin{alignat}{2}
    \dot{\bfK}(t) =\,& \bfF(t, \bfK(t)\bfV_0^{\top})\bfV_0\,,\quad \,&&\bfK(t_0) = \bfU_0\bfS_0\,,\\
    \dot{\bfL}(t) =\,& \bfF(t, \bfU_0\bfL(t)^{\top})^{\top}\bfU_0\,,\quad \,&&\bfL(t_0) = \bfV_0\bfS_0^{\top}\,.
\end{alignat}
The time updated basis matrices at time $t_1 = t_0 + \Delta t$, which we denote by $\wh\bfU$ and $\wh\bfV$ can then be retrieved by a QR decomposition of $[\bfK(t_1), \bfU_0] = \wh\bfU\bfS_K$ and $[\bfL(t_1),\bfV_0] = \wh\bfV\bfS_L$, where $[\bfA, \bfB]\in\mathbb{R}^{m\times 2r}$ denotes the concatenation of two matrices $\bfA$, $\bfB\in\mathbb{R}^{m\times r}$. Lastly, to update the coefficient matrix $\bfS$, the integrator constrains the dynamics to $\{\wh\bfU\bfS \wh\bfV^{\top} \,|\, \bfS \in\R^{2r \times 2r}\}\subset\mathcal{M}_{2r}$ by solving the Galerkin system
\begin{align}
    \dot{\bfS}(t) =\,& \wh\bfU^{\top}\bfF(t, \wh\bfU\bfS(t)\wh\bfV^{\top})\wh\bfV\,,\quad \bfS(t_0) = (\wh\bfU^{\top}\bfU_0)\bfS_0(\bfV_0^{\top}\wh\bfV)\,.
\end{align}
\end{subequations}
The time-updated coefficient matrix $\wh\bfS$ is retrieved from $\bfS(t_1)$. Since the time-updated solution $\wh\bfY = \wh \bfU \wh \bfS \wh \bfV^{\top}$ is of rank $2r$, a new rank can be determined by a truncated singular value decomposition of $\wh\bfS$. Given $\wh\bfS = \bfP\bfSigma\bfQ^{\top}$, where $\bfSigma = \text{diag}(\sigma_1,\cdots,\sigma_{2r})$ collects the singular values of $\wh\bfS$, the singular value decomposition is truncated at a new rank $r_1$ such that $\sum_{i=r_1+1}^{2r} \sigma_i \leq \vartheta$ with a user-determined truncation threshold $\vartheta$, see \cite{ceruti_rank-adaptive_2022}. This yields a new rank $r_1$ factorization $\bfY_1 = \bfU_1\bfS_1\bfV_1^{\top}$. The process of evolving the solution along \eqref{eq:KLS} with subsequent truncation is repeated until the factorized solution at a given final time $t_k$ is reached after $k$ steps.
The equations \eqref{eq:KLS} inherit the stiffness of the original problem, thus, no further stiffness is added from the underlying geometry of $\mathcal{M}$. Therefore, classical time integration methods can be used to evolve the factorized solution. For a detailed discussion of BUG integrators, their inherent parallelism, and extensions to higher order, see \cite{ceruti_unconventional_2022,ceruti_rank-adaptive_2022,ceruti_parallel_2024,ceruti_robust_2024,kusch_second-order_2025}.
\section{Methods}\label{sec:methods}
In this section we detail the numerical method used to solve the problem described in section \ref{sec:2}. We loosely follow the strategy in \cite{kusch_robust_2023} of splitting the problem in collided and uncollided as well as scattering and streaming parts and using a P$_\text{N}$-based dynamical low-rank approach for the collided equations. In contrast to \cite{kusch_robust_2023}, we however use a ray tracer including energy straggling for the uncollided part and the derivation for the collided equations differs due to the material dependency in cross sections and stopping power as well as the higher order spatial discretization. Together with the efficient incorporation of multiple materials as detailed in Section~\ref{sec:multiple_materials}, this enables the derivation of an efficient and accurate numerical solver to simulate proton radiation in a realistic setting. 

\subsection{Collided-uncollided split}
We start by splitting equation \eqref{eq:finalContinuous} into a collided and uncollided part using a first collision source method. This will facilitate the inclusion of straggling and avoid challenging boundary conditions in the DLRA part of the method. 
We now separate the angular flux into collided particles $\psi_c$ and uncollided particles $\psi_u$ such that $\psi = \psi_u + \psi_c$. For collided particles, we assume that straggling can be neglected, that is,
\begin{align}
    \frac{\partial^2 T\psi_c}{\partial E^2} \approx 0.
\end{align}
Then, we obtain 
\begin{subequations}
    \begin{align}
        \boldsymbol\Omega  \cdot\nabla \psi_u  =\,& \frac{\partial \mathcal{S}\psi_u}{\partial E} + \frac12 \frac{\partial^2 T\psi_u}{\partial E^2} +  \sum_{i\in \mathbb{A}} w_i \Gamma_{\mathrm{out};i}\psi_u(E,\mathbf{r},\mathbf{\Omega}) , \label{eq:final_uncollided}\\
        -\partial_E \mathbf{\widetilde{\psi}_c}(E,\mathbf{r}) =\,&-\boldsymbol\Omega \cdot\nabla \frac{\mathbf{\widetilde{\psi}_c}(E,\mathbf{r})}{\mathcal{S}(E,\mathbf{r})} + 
 \sum_{i\in \mathbb{A}} w_i \Gamma_i \frac{\mathbf{\widetilde{\psi}_c}(E,\mathbf{r})}{\mathcal{S}(E,\mathbf{r})} +  \sum_{i\in \mathbb{A}} w_i \Gamma_{\mathrm{in};i}  \psi_u(E,\mathbf{r},\mathbf{\Omega}),\;\label{eq:final_collided_multmat}
    \end{align}
\end{subequations}
where $ \Gamma_{\mathrm{out};i}$ and $\Gamma_{\mathrm{out};i}$ are the in- and outscattering parts of the collision operator, respectively. So,
\begin{align}
    \Gamma_{\mathrm{in};i}\psi &=  \int_{\mathbb{S}^2} \Sigma_{s;i}(E,\mathbf{r},\mathbf{\Omega}'\cdot \mathbf{\Omega})\psi(E,\mathbf{r},\mathbf{\Omega}')\mathrm{d}\mathbf{\Omega}'
 , \\ 
  \Gamma_{\mathrm{out};i}\psi &=  \Sigma_{t;i}(E,\mathbf{r})\psi(E,\mathbf{r},\mathbf{\Omega})
\end{align}
for the Boltzmann operator and 
\begin{align}
    \Gamma_{\mathrm{in};i}\psi =  \Gamma_{i}\psi(E,\mathbf{r},\mathbf{\Omega}) ,
    \qquad \Gamma_{\mathrm{out};i}\psi = 0
\end{align}
for the Fokker-Planck operator.
\subsection{Uncollided part: Ray tracer}
For the uncollided particles, we use an in-house analytical ray tracer implemented in fortran \cite[see also][]{burlacu_deterministic_2023,lathouwers_deterministic_2023}. 

Having split off the singular part as described above, the PDE along the ray for $\psi_u$ reads
\begin{equation}
  \frac{\partial \psi_u(E,z)}{\partial z} + \sum_{i\in \mathbb{A}} w_i(z) \Gamma_{\mathrm{out};i}\psi_u(E,z) = \frac{\partial \mathcal{S}^* (E,z) \psi_u(E,z) }{\partial E} + 
  \frac{1}{2} \frac{\partial}{\partial E} T(E,z) \frac{\partial \psi_u(E,z)}{\partial E},
\end{equation}
where $z$ denotes the one-dimensional spatial position, i.e. depth along the ray, and we define $\mathcal{S}^* = \mathcal{S} + \frac{1}{2} \frac{dT}{dE}$ to obtain the standard diffusion form of straggling.
We use a discontinuous Galerkin scheme to discretize the PDE in energy. The energy domain is divided into 128 equally-sized groups. 
The slowing down term is based on the Lax-Friedrichs scheme, whereas the straggling term utilizes the SIPG method. Polynomials up to
second order are used as basis vectors (3 dof per energy element).
After discretization in $E$, the system can be written as
\begin{equation}
 \mathbf M_E \frac{\partial \boldsymbol \psi_u(z)}{\partial z} + \mathbf G \boldsymbol \psi_u(z) = 0,
\end{equation}
where $ \mathbf M_E$ is the mass matrix for the energy basis and $ \mathbf  G$ contains elements corresponding to absorption, slowing down and straggling.
Here, $\psi_u(z=0) = f(E)$ is the projection of a given energy distribution (Gaussian).
We use a Crank-Nicolson scheme to solve this ODE system along the depth. The step size is determined such that it is smaller than the maximum
step length (0.01 cm). The unscattered flux distribution is computed at the cell midpoints of the mesh. Based on the spatial and angular distribution of the beam source, a collection of start and end points for rays is generated. The solution is then summed over all ray positions within a beam source and the resulting space and energy dependent flux acts as a source term for the collided equation through inscattering.

\subsection{Collided part: Dynamical low-rank approximation}
We now want to use DLRA to efficiently solve the costly collided part of the equation. Since DLRA is a method designed for time-dependent problems, we now transform from energy $E$ to pseudo-time $t(E)=E_{max}-E$, resulting in the following equation:
\begin{align}
 \partial_t \mathbf{\widetilde{\psi}_c}(t,\mathbf{r}) =\,&-\boldsymbol\Omega \cdot\nabla \frac{\mathbf{\widetilde{\psi}_c}(t,\mathbf{r})}{\mathcal{S}(t,\mathbf{r})} + 
 \sum_{i\in \mathbb{A}} w_i \Gamma_i \frac{\mathbf{\widetilde{\psi}_c}(t,\mathbf{r})}{\mathcal{S}(t,\mathbf{r})} +  \sum_{i\in \mathbb{A}} w_i \Gamma_{\mathrm{in};i}  \psi_u(t,\mathbf{r},\mathbf{\Omega}).\;
\end{align}
Next, we will discretize this equation in angle and space to obtain a (pseudo)time dependent ordinary matrix differential equation.

\subsubsection{Angular discretization}
\label{sec:Pn}
We use a modal discretization of the angular variable. To this end, we define the real-valued spherical harmonics basis 
\begin{align*}
    m_{\ell}^k = 
    \begin{cases}
        \frac{(-1)^k}{\sqrt{2}}\left( Y_{\ell}^k + (-1)^k Y_{\ell}^{-k} \right), & k > 0\;, \\
        Y_{\ell}^0 & k = 0 \;, \\
        -\frac{(-1)^k i}{\sqrt{2}}\left( Y_{\ell}^{-k} - (-1)^k Y_{\ell}^{k} \right), & k < 0\;,
    \end{cases}
\end{align*}
which uses the spherical harmonics basis functions
\begin{align*}
Y_{\ell}^k(\boldsymbol\Omega) = \sqrt{\frac{2\ell +1}{4\pi}\frac{(\ell-k)!}{(\ell+k)!}}\ e^{ik\varphi}P_{\ell}^k(\mu)\;,
\end{align*}
where $P_{\ell}^k$ are the associated Legendre polynomials. Collecting all basis functions up to degree $N$ in a vector
\begin{align*}
    \mathbf m = (m_0^0, m_1^{-1}, m_1^{0}, m_1^{1},\cdots, m_N^{N})^{\top}\in\mathbb{R}^{m}, \; m= (N+1)^2
\end{align*}
gives the modal approximation $\psi_c(t,\mathbf{r},\boldsymbol\Omega) \approx \mathbf u(t,\mathbf{r})^{\top}\mathbf m(\boldsymbol\Omega)$. Since we deal with monodirectional beam sources in this work, for a specific source with beam direction $\boldsymbol{\Omega}_{in}$, the uncollided flux $\psi_u (t,\mathbf{r},\boldsymbol{\Omega})= \psi_u(t,\mathbf{r})\delta(\boldsymbol\Omega - \boldsymbol\Omega_{\mathrm{in}})$ at a specific point in space and time/energy is a scalar and the transformation from nodal to modal representation in angle can be done using the vector $\mathbf T_M = \mathbf{m}(\boldsymbol{\Omega}_{in}) $. For a generalization to beams with a distribution in angle see \cite{kusch_robust_2023}. Then, the P$_N$ equations read
\begin{align}
\label{eq:discAngle}
    \partial_t \mathbf u (t,\mathbf{r}) &=-\mathbf A\cdot\nabla \frac{\mathbf u(t,\mathbf{r})}{\mathcal{S}(t,\mathbf{r})}- \sum_{i\in \mathbb{A}}w_i(\mathbf r)\Sigma_{t;i}(t)\frac{\mathbf u (t,\mathbf{r})}{\mathcal{S}(t,\mathbf{r})}\nonumber \\
    &+  \sum_{i\in \mathbb{A}}w_i(\mathbf r)\mathbf G_i(t)\frac{\mathbf u(t,\mathbf{r})}{\mathcal{S}(t,\mathbf{r})} +  \sum_{i\in \mathbb{A}}w_i(\mathbf r)\mathbf G_i(t)\mathbf T_M\psi_u (t,\mathbf{r}),\;
\end{align}
where $\mathbf A\cdot\nabla := \mathbf A_1\partial_{x} + \mathbf A_2\partial_y+ \mathbf A_3\partial_z$ with $\mathbf A_i := \int_{\mathbb{S}^2} \mathbf m \mathbf m^{\top} \Omega_i \, d\boldsymbol\Omega$.

The diagonal scattering matrix of material $i$ is defined as $\mathbf G_i(t)=\left(G_{ i; pq}\right)_{1\leq p,q\leq m}$ where $$G_{ i; pp}(t, \mathbf r) = 2\pi\int_{[-1,1]}P_{k}(\mu_0)\Sigma_{s;i}(t,\mu_0)\,d\mu_0,$$ for Boltzmann and $$G_{ i; pp}(t) = -\frac{\xi_{1;i}}{2} k\cdot(k+1),$$ for Fokker-Planck, with $p = l^2 + l + k, \; 0\leq k \leq N, \; |l|\leq k$. The total cross section is then $\Sigma_{t;i}(t) = G_{i;11}(t) $. 
Note, that we further apply the extended transport correction to the scattering moments in the Boltzmann equation \cite{drumm_analysis_2007} and the correction described in \cite{landesman_angular_1989,morel_fokker-planck_1981} for the Fokker-Planck moments. Both methods aim at a "delta-function corrected expansion" \cite{landesman_angular_1989}, where the extended transport correction is chosen such that the truncated expansion matches moments 0 to $N+1$ exactly \cite{drumm_analysis_2007}.
\subsubsection{Spatial discretization}
In this work, we consider a full three-dimensional spatial domain, i.e., the spatial variable is $\mathbf{r}=(x,y,z)^T \in \mathbb{R}^3$. In order to sufficiently resolve the sharp gradients of the energy deposition near the Bragg peak, a higher-order method is necessary \cite{lathouwers_angular_2019}. Here, we use a second order upwind scheme. For this, the spatial domain is split into $n=n_x\cdot n_y \cdot n_z$ cells in the form of a structured rectangular grid with uniform one-dimensional grids $x_1\leq x_2 \leq ... \leq x_{n_x}$ , $y_1\leq y_2 \leq ... \leq y_{n_y}$  and $z_1\leq  z_2 \leq ... \leq z_{n_z}$ with grid sizes $\Delta_x, \Delta_y,\Delta_z$ respectively. 

The cell of index $\mathrm{idx}(i,j,k):=(k-1)\cdot n_x\cdot n_y+(j-1)\cdot n_x+i$ is defined on $I_{i,j,k}:= [x_i,x_{i+1}]\times[y_i,y_{i+1}]\times [z_i,z_{i+1}]$ with $ i=1,...,n_x, \; j=1,...,n_y, \; k=1, ...,n_z$. Assuming we want to solve the semi-discrete equation (\ref{eq:discAngle}) from the previous section, we now collect the numerical solution in a matrix $\mathbf{u}(t)\in\mathbb{R}^{n\times m}$ with entries $u_{idx(i,j,k),p}(t), \; p=1,...,m$. We define the second order upwind difference stencils $\mathbf{D}_{x,y,z}^{+,-} \in \mathbb{R}^{n\times n}$ as
\begin{align*}
    \mathbf{D}_{x;\mathrm{idx}(i,j,k),\mathrm{idx}(i,j,k)}^{+,-} &= \frac{\pm 3}{2\Delta_x}, \; \mathbf{D}_{x;\mathrm{idx}(i,j,k),\mathrm{idx}(i\pm 1,j,k)}^{+,-} = \frac{\mp 4}{2\Delta_x},  \;  \mathbf{D}_{x;\mathrm{idx}(i,j,k),\mathrm{idx}(i\pm 2,j,k)}^{+,-} = \frac{\pm 1}{2\Delta_x} \\
    \mathbf{D}_{y;\mathrm{idx}(i,j,k),\mathrm{idx}(i,j,k)}^{+,-} &= \frac{\pm 3}{2\Delta_y}, \;  \mathbf{D}_{y;\mathrm{idx}(i,j,k),\mathrm{idx}(i,j\pm 1,k)}^{+,-} = \frac{\mp 4}{2\Delta_y},   \; \mathbf{D}_{y;\mathrm{idx}(i,j,k),\mathrm{idx}(i,j\pm 2,k)}^{+,-} = \frac{\pm 1}{2\Delta_y} \\
    \mathbf{D}_{z;\mathrm{idx}(i,j,k),\mathrm{idx}(i,j,k)}^{+,-} &= \frac{\pm 3}{2\Delta_z}, \;   \mathbf{D}_{z;\mathrm{idx}(i,j,k),\mathrm{idx}(i,j,k\pm1)}^{+,-} = \frac{\mp 4}{2\Delta_z},   \;  \mathbf{D}_{z;\mathrm{idx}(i,j,k),\mathrm{idx}(i,j,k\pm2)}^{+,-} = \frac{\pm 1}{2\Delta_z}
\end{align*}

Further, we collect the eigenvectors of $\mathbf{A}_{x,y,z} = \mathbf{V}_{x,y,z}\mathbf{\Lambda}_{x,y,z} \mathbf{V}_{x,y,z}^T$ in $\mathbf{V}_{x,y,z}$ and split the eigenvalues collected in $\mathbf{\Lambda}_{x,y,z} = \mathrm{diag}(\lambda_1^{x,y,z},...,\lambda_m^{x,y,z})$ into $\left(\mathbf{\Lambda}_{x,y,z;p,p}^+\right)_{p=1,...,m}=\mathrm{max}(0,\lambda_p^{x,y,z})$ and $\left(\mathbf{\Lambda}_{x,y,z;p,p}^-\right)_{p=1,...,m}=\mathrm{min}(0,\lambda_p^{x,y,z})$. Then we obtain a large matrix differential equation of the form 
\begin{align}
\label{eq:streamscattsplit}
   \mathbf{\dot{u}}(t) = \mathbf{F_S}(\mathbf{u}(t)) + \mathbf{G}(t,\mathbf{u}(t)), 
\end{align}
where 
\begin{align*}
    \mathbf{F_S}(\mathbf{u}(t)) &:= \left(\mathbf D_x^+\left(\mathbf V_x \boldsymbol{\mathcal{S}}^{-1}(t) \mathbf{u}(t)\right)\mathbf\Lambda_x^+ + \mathbf D_x^-\left( \mathbf V_x \boldsymbol{\mathcal{S}}^{-1}(t) \mathbf{u}(t)\right) \mathbf \Lambda_x^-\right) \mathbf V_x^T \\ &+ \left(\mathbf D_y^+\left(\mathbf V_y\boldsymbol{\mathcal{S}}^{-1}(t) \mathbf{u}(t)\right) \mathbf \Lambda_y^+ + \mathbf D_y^-\left(\mathbf V_y \boldsymbol{\mathcal{S}}^{-1}(t)\mathbf{u}(t)\right) \mathbf \Lambda_y^-\right) \mathbf V_y^T\\ &+\left(\mathbf D_z^+\left(\mathbf V_z \boldsymbol{\mathcal{S}}^{-1}(t)\mathbf{u}(t)\right)\mathbf \Lambda_z^+ +\mathbf D_z^-\left(\mathbf V_z \boldsymbol{\mathcal{S}}^{-1}(t) \mathbf{u}(t)\right) \mathbf \Lambda_z^-\right) \mathbf V_z^T \\
    \mathbf{G}(t,\mathbf u(t)) &:=  \sum_{i\in \mathbb{A}} \mathbf w_i \boldsymbol{\mathcal{S}}^{-1}(t) \mathbf{u}(t) \mathbf\Sigma_{t;i} +  \sum_{i\in \mathbb{A}} \mathbf w_i \boldsymbol{\mathcal{S}}^{-1}(t)  \mathbf{u}(t) \mathbf G_i\\ &+  \sum_{i\in \mathbb{A}} \mathbf w_i  \boldsymbol{\psi}_u(t) \mathbf{T_M} \mathbf G_i,
\end{align*}
here $\mathbf w_i \text{ and } \boldsymbol{\mathcal{S}}(t) \in \mathbb{R}_+^{n\times n}$ are diagonal matrices containing the material composition weights and stopping power at each spatial grid point, respectively. Here, $\boldsymbol{\psi}_u$ is the uncollided solution from the ray tracer evaluated at the (same) spatial grid points.
Since solving this at sufficiently high resolution is extremely costly we use the dynamical low-rank approximation to evolve the solution in (pseudo-)time.
\subsubsection{Energy/time update using the dynamical low-rank approximation}
For this, we further split equation \eqref{eq:streamscattsplit} into its streaming and scattering parts using a first-order Lie splitting, i.e.,
\begin{subequations}
    \begin{alignat}{2}
    \mathbf{\dot{u}}_{\Romannum{1}}(t) =\,& \mathbf{F_S}(\mathbf{u}_{\Romannum{1}}(t))\,\quad &&\mathbf{u}_{\Romannum{1}}(t_0) = \mathbf{u}(t_0)\,, \label{eq:split-streaming}\\
    \mathbf{\dot{u}}_{\Romannum{2}}(t) =\,&  \mathbf{G}(t,\mathbf{u}_{\Romannum{2}}(t))\,\quad &&\mathbf{u}_{\Romannum{2}}(t_0) = \mathbf{u}_{\Romannum{1}}(t_1)\,.\label{eq:split-scattering}
\end{alignat}
\end{subequations}
 For the streaming update \eqref{eq:split-streaming} we use the augmented BUG integrator from section \ref{sec:DLRA_intro} to update $\mathbf{u}_{\Romannum{1}}$ in factorized form from $\mathbf U^0,\mathbf S^0,\mathbf V^0$ to $\mathbf U^{\frac{1}{2}},\mathbf S^{\frac{1}{2}},\mathbf V^{\frac{1}{2}}$, by simply plugging in $\mathbf F_S$ for the right hand side and using a fourth-order Runge-Kutta method for the time updates. Note, that we still do a full time step from $t_0$ to $t_1=t_0 + \Delta t$ and the $\frac{1}{2}$ in the exponent merely denotes that we have only completed the streaming part of the update, i.e., $\mathbf U^{\frac{1}{2}}\mathbf S^{\frac{1}{2}}\mathbf V^{\frac{1}{2},\top}\approx \mathbf{u}_{\Romannum{1}}(t_1)$. The fourth-order Runge--Kutta method is primarily used to stabilize the stiff streaming operators in $\mathbf F_S$.

   For the scattering update, we use the matrix projector splitting integrator for the collided flux and the augmented BUG integrator on the remainder according to \cite{kusch_stability_2023}. This reduces the update to a K-, S- and L-step for inscattering from the uncollided flux and just the L-step for self-scattering of the collided flux. Thus, to perform the scattering update \eqref{eq:split-scattering}, we need to solve the following differential equations
    \begin{enumerate}
    \item \textbf{$L$-step, collided:} \begin{align*}
        \dot{\mathbf L}(t) &= -  \sum_{i\in \mathbb{A}}  \mathbf{U}^{\frac{1}{2},\top} \mathbf w_i  \boldsymbol{\mathcal{S}}^{-1}(t)\mathbf{U}^{\frac{1}{2}} \Sigma_{t,i}\mathbf{L}(t) +  \sum_{i\in \mathbb{A}} \mathbf{U}^{\frac{1}{2},\top}\mathbf w_i\boldsymbol{\mathcal{S}}^{-1}(t)\mathbf{U}^{\frac{1}{2}} \mathbf{L}(t)\mathbf{G}_i.\\ \mathbf L(t_0) &= \mathbf{S}^{\frac{1}{2}} \mathbf{V}^{{\frac{1}{2}},\top}.
    \end{align*}
    Determine the updated basis using a QR decomposition of $\mathbf{L}(t_1) =\widetilde{\mathbf{S}}^{\frac{1}{2}}  \widetilde{\mathbf{V}}^{{\frac{1}{2}}}$.
		\item \textbf{$K$-step, uncollided}: 
			\begin{align*}
			\dot{\mathbf K}(t) &=   \sum_{i\in \mathbb{A}} \mathbf w_i \boldsymbol{\mathcal{S}}^{-1}\boldsymbol{\psi}_u  \mathbf T_M\mathbf{G}_i \mathbf{V}^{\frac{1}{2}}, \\
            \mathbf K(t_0) &= \mathbf{U}^{\frac{1}{2}}\mathbf{S}^{\frac{1}{2}}
			\end{align*}
			Determine the updated and augmented basis at $t_1=t_0 + \Delta t$ using a QR decomposition of $[\bfK(t_1), \bfU^\frac{1}{2}] = \wh\bfU\bfS_K$.
		\item \textbf{$L$-step, uncollided}: 
			\begin{align*}
			\dot{\mathbf L}(t) &=  \sum_{i\in \mathbb{A}} \mathbf{U}^{{\frac{1}{2}},\top} \mathbf w_i \boldsymbol{\mathcal{S}}^{-1}\boldsymbol{\psi}_u \mathbf T_M \mathbf{G}_i, \\ \mathbf L(t_0) &= \widetilde{\mathbf{S}}^{\frac{1}{2}}  \widetilde{\mathbf{V}}^{{\frac{1}{2}}}.
			\end{align*}
			Determine the updated and augmented basis at $t_1=t_0 + \Delta t$ using a QR decomposition of $[\bfL(t_1),\widetilde{\mathbf{V}}^{{\frac{1}{2}}}] = \wh\bfV\bfS_L$.
		\item \textbf{$S$-step, uncollided}:
			\begin{align*}
			\dot{\mathbf S}(t) &=  \sum_{i\in \mathbb{A}}  \wh\bfU^{\top} \mathbf w_i \boldsymbol{\mathcal{S}}^{-1}\boldsymbol{\psi}_u \mathbf T_M \mathbf{G}_i\wh\bfV, \\ \quad \bfS(t_0) &= (\wh\bfU^{\top}\bfU^{\frac{1}{2}})\bfS^{\frac{1}{2}}(\widetilde{\mathbf{V}}^{{\frac{1}{2}},\top}\wh\bfV)
			\end{align*}
			and set $\mathbf S^1 = \mathbf S(t_1)$.
	\end{enumerate}
    Similar to \cite{kusch_robust_2023}, we use an implicit Euler method for the first step to treat stiff scattering terms and an explicit Euler method for the time updates in the remaining steps. 

\section{Numerical results}
\label{sec:numResults}
In the following, we present a number of numerical experiments in homogeneous and heterogeneous media to validate and compare our method against an otherwise equivalent full-rank method as well as a state-of-the-art Monte Carlo code (TOPAS MC) \cite{perl_topas_2012} with only electromagnetic interactions. We will first consider a single uni-directional beam in different media and later add a second beam with a different incoming angle to investigate the effect this has on the chosen ranks. All beams have a Gaussian spatial distribution in the two-dimensional plane perpendicular to the beam direction with 3mm standard deviation and the intial beam energy follows a Gaussian distribution with standard deviation $1\%$ of the mean beam energy.
\begin{figure}[h!]
    \centering
    \begin{subfigure}{0.35\linewidth}
        \includegraphics[width=\linewidth]{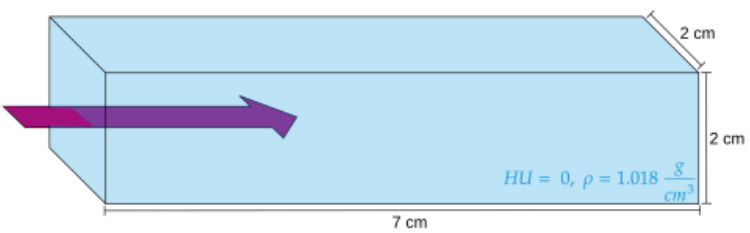}
        \caption{}
           \label{fig:testCasesA}
    \end{subfigure}
   \begin{subfigure}{0.35\linewidth}
        \includegraphics[width=\linewidth]{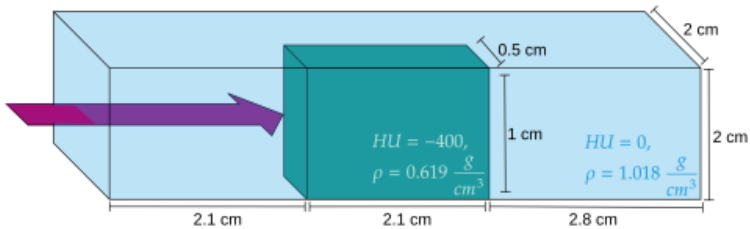}
        \caption{}
           \label{fig:testCasesB}
    \end{subfigure}
    \begin{subfigure}{0.25\linewidth}
        \includegraphics[width=\linewidth]{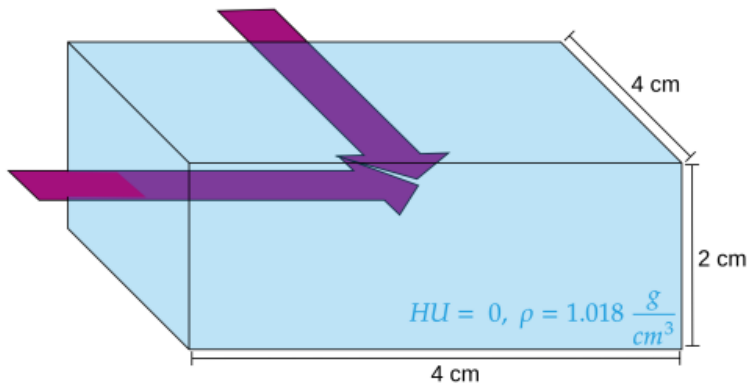}
        \caption{}
           \label{fig:testCasesC}
    \end{subfigure}
    \caption{Schematic representations of the domain size, composition and beam directions in the three numerical test cases: (a) homogeneous, water-like $2\times2\times7$ cm domain with single beam in z-direction hitting middle of x-y-plane (b) heterogeneous $2\times2\times7$ cm domain representing water with a lung insert with a single beam identical to (a) and (c) homogeneous, water-like $2\times4\times4$ cm domain with two perpendicular beams pointing in z- and negative y-direction. }
    \label{fig:testCases}
\end{figure}
Since the problem at full angular and spatial resolution is infeasible to solve at full rank, we will validate the low-rank method against the deterministic full-rank method on a coarser grid and then compare a higher resolution low-rank solution to the Monte Carlo reference. For the spatial resolution we always use a coarser gridsize of 1 mm and a finer one with 0.25 mm. Since initial numerical experiments indicate that more angular modes are required when using the Boltzmann collision operator, we choose $N=75$ for the coarser and $N=115$ for the finer resolution here. In case of the Fokker-Planck approximation, which uses a simplified collision operator, $N=19$ for the coarser and $N=75$ for the finer resolution have proven sufficient. The truncation tolerance used to adapt the rank is chosen as $\vartheta=0.01$ in all cases.

\subsection{Homogeneous Case}
\begin{figure}[h!]
    \centering
    \begin{subfigure}{0.225\linewidth}
    \centering
        \includegraphics[width=\linewidth]{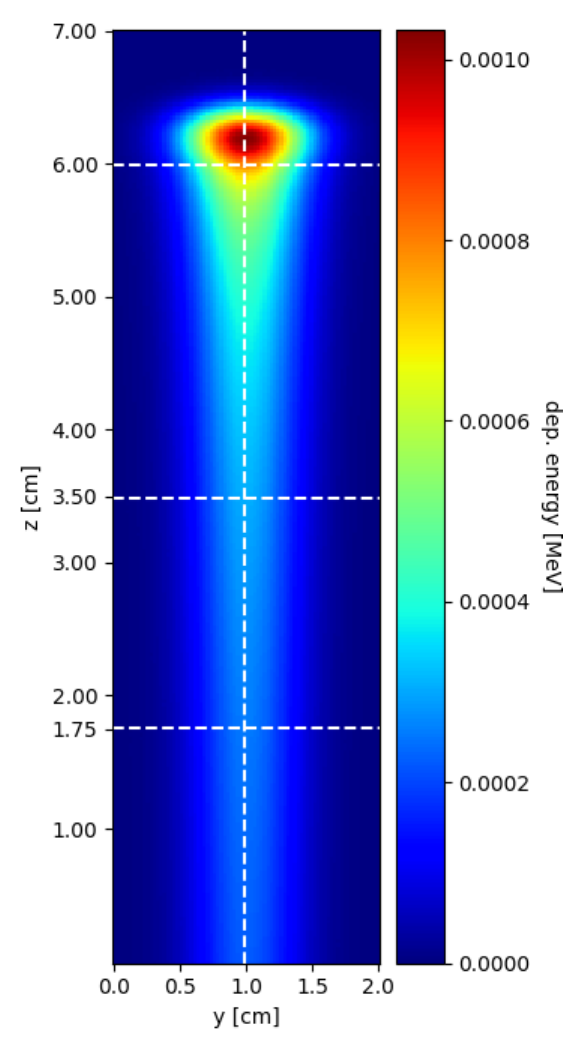}
        \caption{Monte Carlo, \\1e8 histories}
    \end{subfigure}
   \begin{subfigure}{0.225\linewidth}
   \centering
        \includegraphics[width=\linewidth]{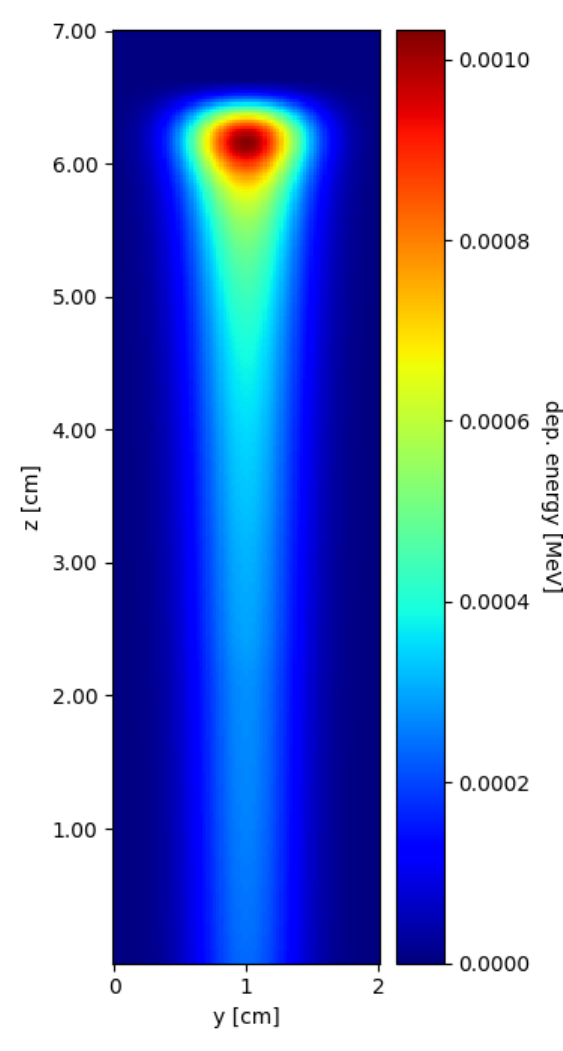}
        \caption{DLRA, P$_\text{115}$,\\  $\Delta_x=0.25$ mm}
    \end{subfigure}
    \begin{subfigure}{0.225\linewidth}
    \centering
        \includegraphics[width=\linewidth]{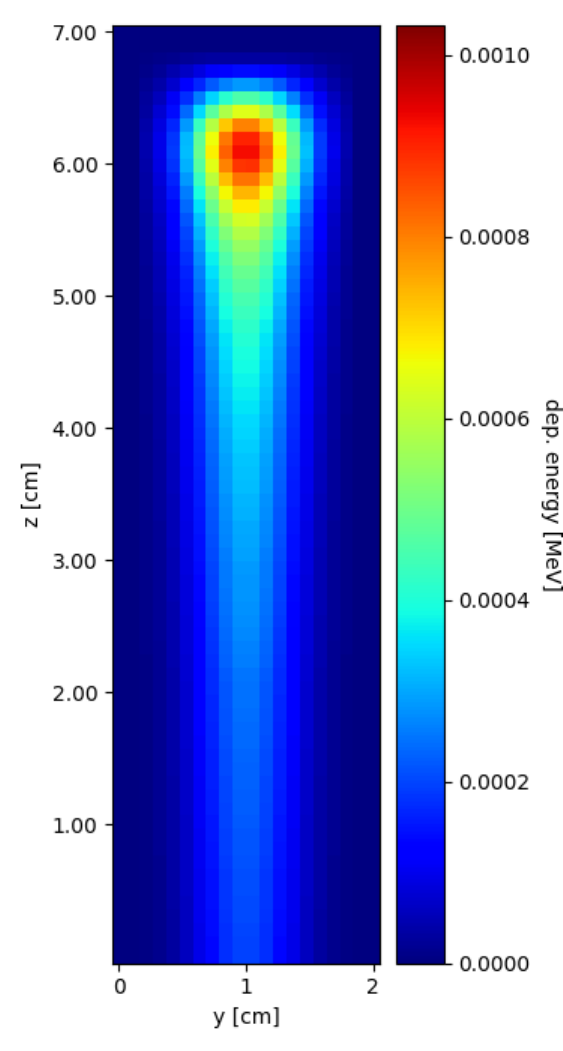}
        \caption{DLRA, P$_\text{75}$,\\ $\Delta_x=1$ mm}
    \end{subfigure}
    \begin{subfigure}{0.225\linewidth}
    \centering
        \includegraphics[width=\linewidth]{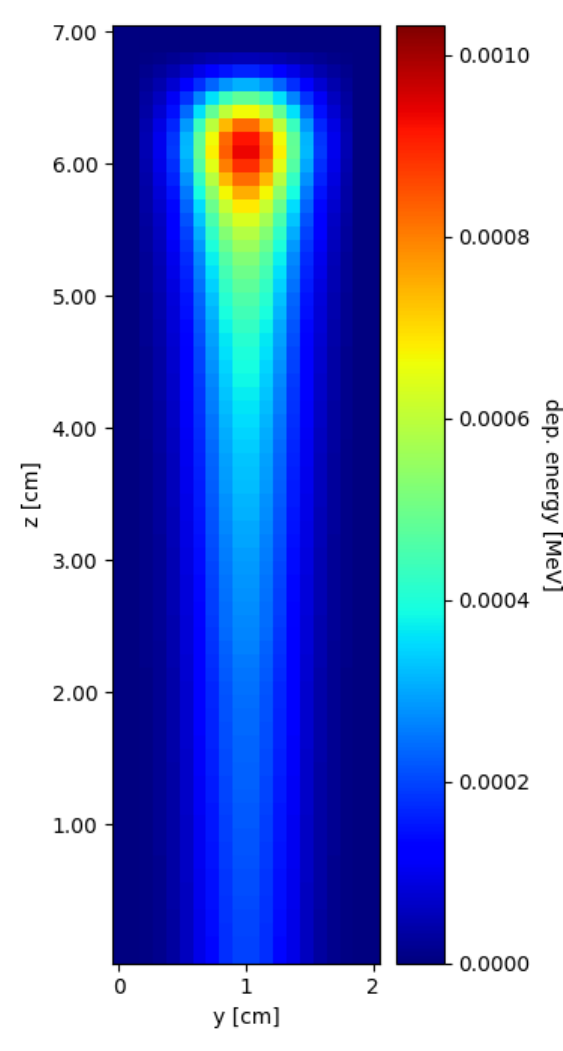}
        \caption{full rank P$_\text{75}$,\\  $\Delta_x=1$ mm}
    \end{subfigure}\\
        \begin{subfigure}{0.31\linewidth}
        \centering
        \includegraphics[width=\linewidth]{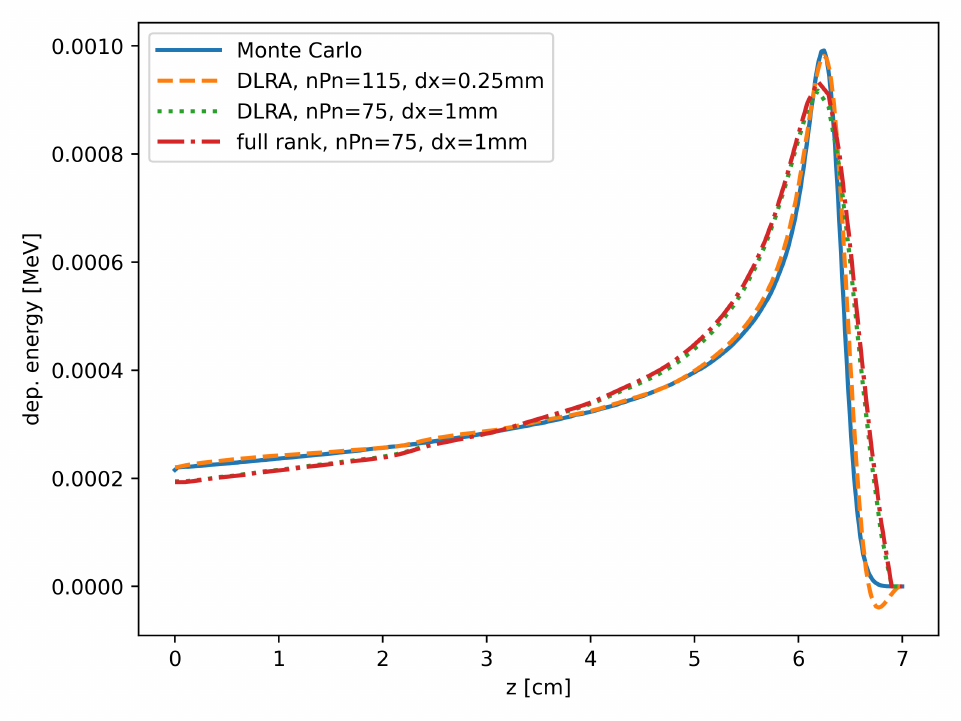}
        \caption{Depth.}
    \end{subfigure}
   \begin{subfigure}{0.31\linewidth}
   \centering
        \includegraphics[width=\linewidth]{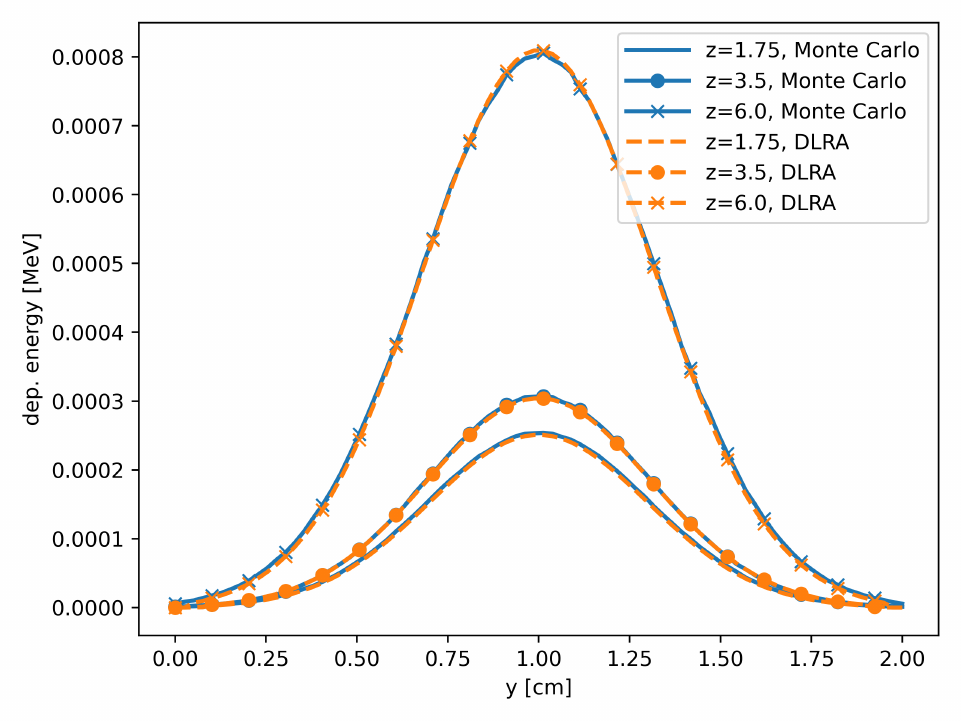}
        \caption{Lateral, MC vs DLRA.}
    \end{subfigure}
 \begin{subfigure}{0.31\linewidth}
 \centering
        \includegraphics[width=\linewidth]{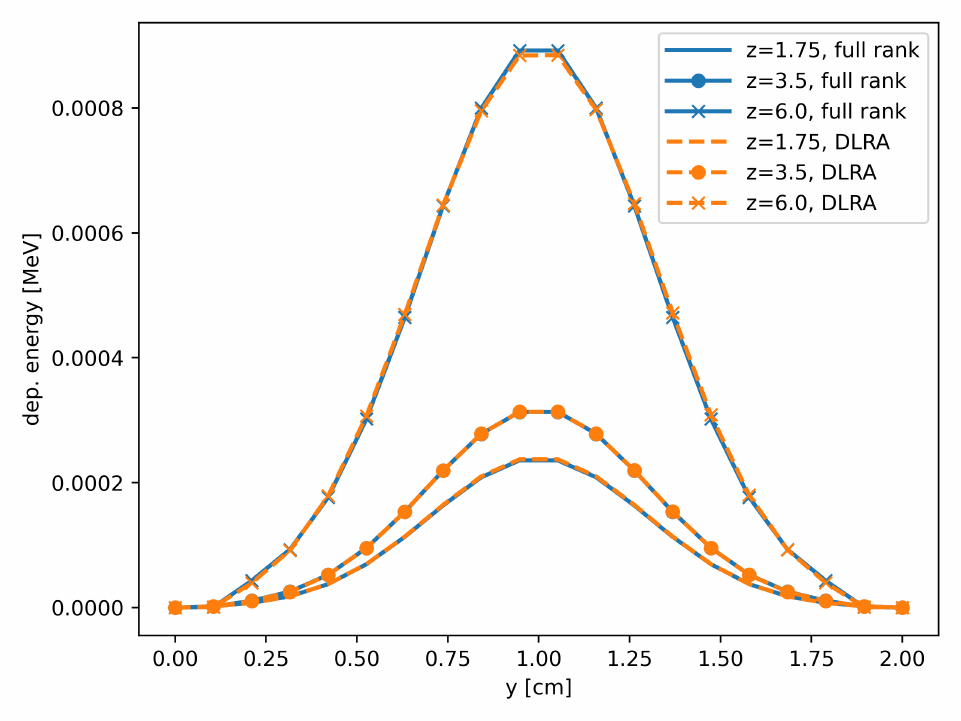}
        \caption{Lateral, DLRA vs full rank.}
    \end{subfigure}
    \caption{Surface plots (a)-(d) and profiles (e)-(g) of deposited energy in the homogeneous test case with a single beam for the Monte Carlo reference  compared to DLRA at high and low spatial and angular resolution and the full-rank deterministic method based on the Boltzmann equation. Lines in (a) indicate positions of profiles plotted in (e)-(g).}
    \label{fig:surfPlotsHomBoltzmann}
\end{figure}
First, we consider a 90 MeV proton beam in a homogeneous $2\times2\times7$ cm box (see figure \ref{fig:testCasesA}). The waterlike material has 0 Hounsfield Units (HU). These values are converted to material composition weights $\mathbf w_i$ and a density of 1.018 $\frac{\text{g}}{\text{cm}^3}$ using the conversion described in \cite{schneider_correlation_2000}.
   \begin{figure}[h!] 
     \centering
     \begin{subfigure}[b]{0.425\textwidth}
         \centering
         \includegraphics[width=\textwidth]{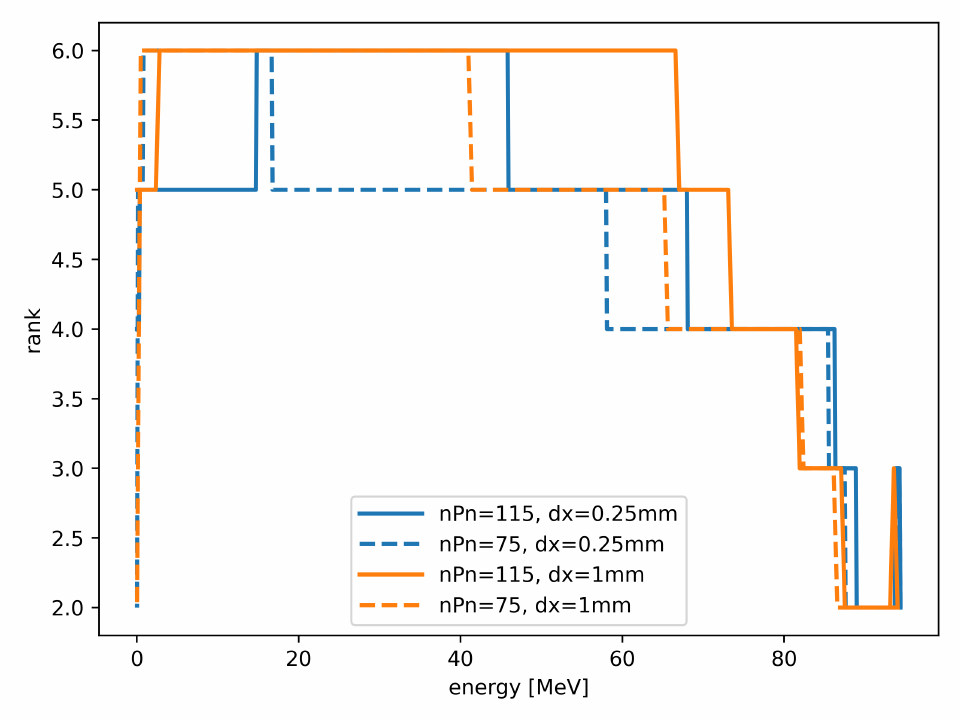}
         \caption{Homogeneous test case.}
     \end{subfigure}
          \begin{subfigure}[b]{0.425\textwidth}
         \centering
         \includegraphics[width=\textwidth]{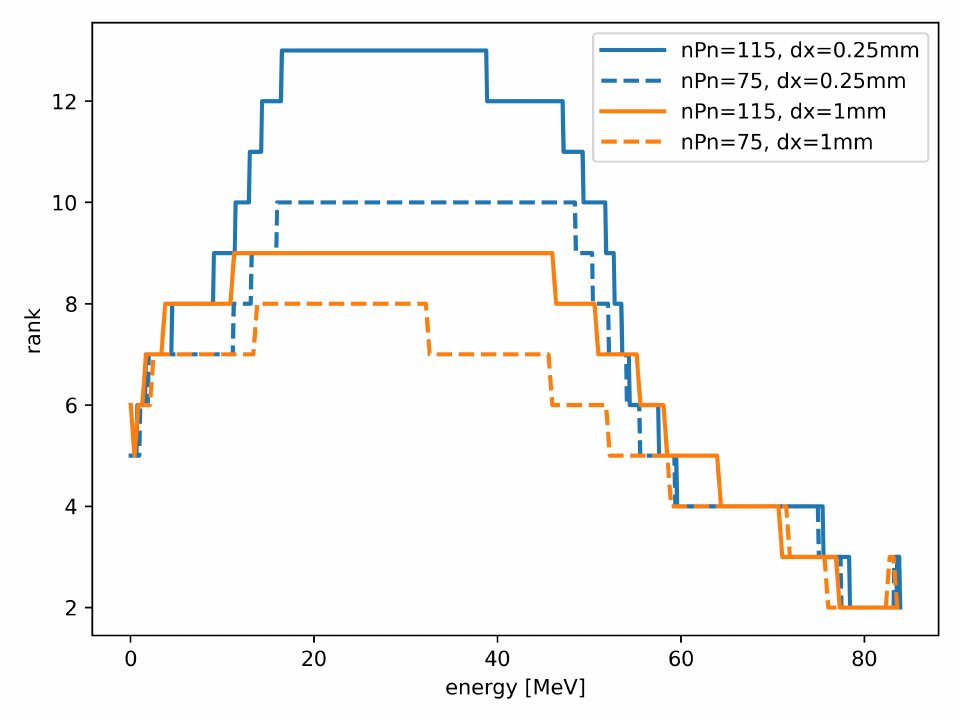}
         \caption{Heterogeneous test case.}
     \end{subfigure}
        \caption{Ranks chosen by rank-adaptive integrator for different spatial and angular resolutions and the Boltzmann equation in (a) the homogeneous and (b) heterogeneous test case with a single beam.}
        \label{fig:RanksBoltzmann}
\end{figure} 
\begin{figure}[h!]
    \centering
    \begin{subfigure}{0.225\linewidth}
     \centering
        \includegraphics[width=\linewidth]{Fig3a_5a.pdf}
        \caption{Monte Carlo, \\1e8 histories}
    \end{subfigure}
   \begin{subfigure}{0.225\linewidth}
    \centering
        \includegraphics[width=\linewidth]{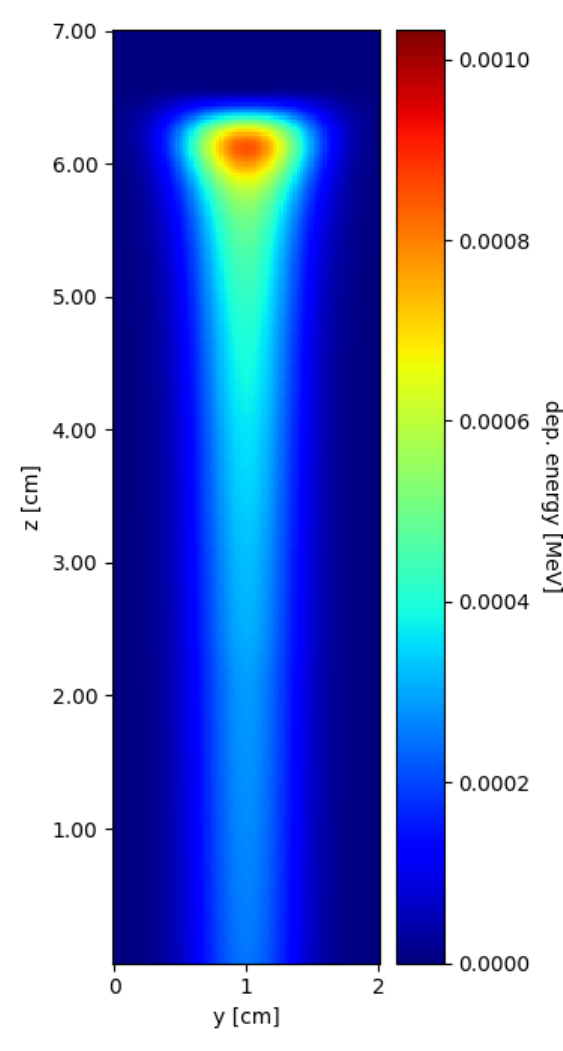}
        \caption{DLRA, P$_\text{75}$, \\ $\Delta_x=0.25$ mm}
    \end{subfigure}
    \begin{subfigure}{0.225\linewidth}
     \centering
        \includegraphics[width=\linewidth]{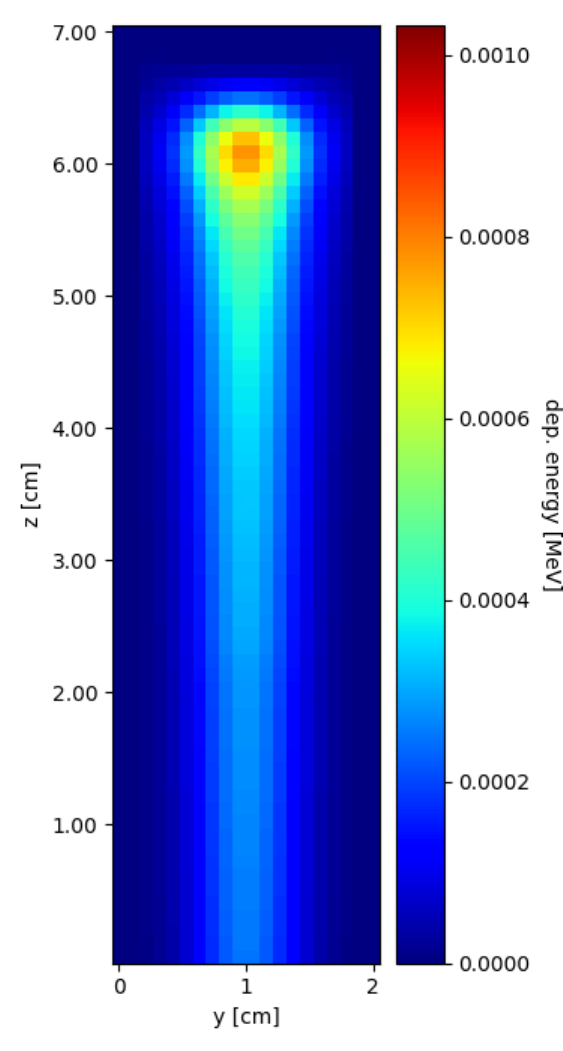}
        \caption{DLRA, P$_\text{19}$,\\  $\Delta_x=1$ mm}
    \end{subfigure}
    \begin{subfigure}{0.225\linewidth}
     \centering
        \includegraphics[width=\linewidth]{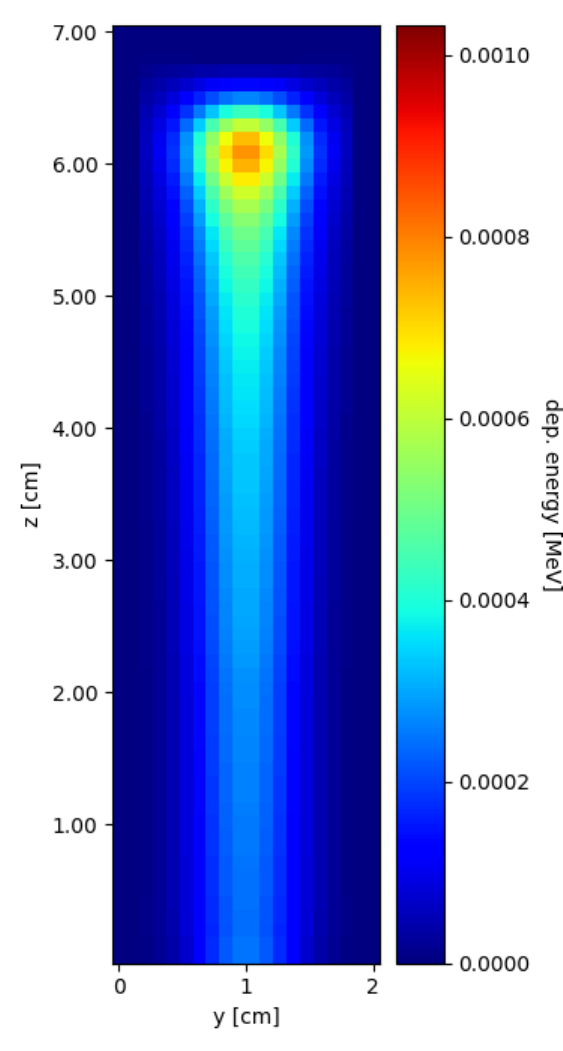}
        \caption{full rank P$_\text{19}$, \\ $\Delta_x=1$ mm}
    \end{subfigure} \\
            \begin{subfigure}{0.31\linewidth}
             \centering
        \includegraphics[width=\linewidth]{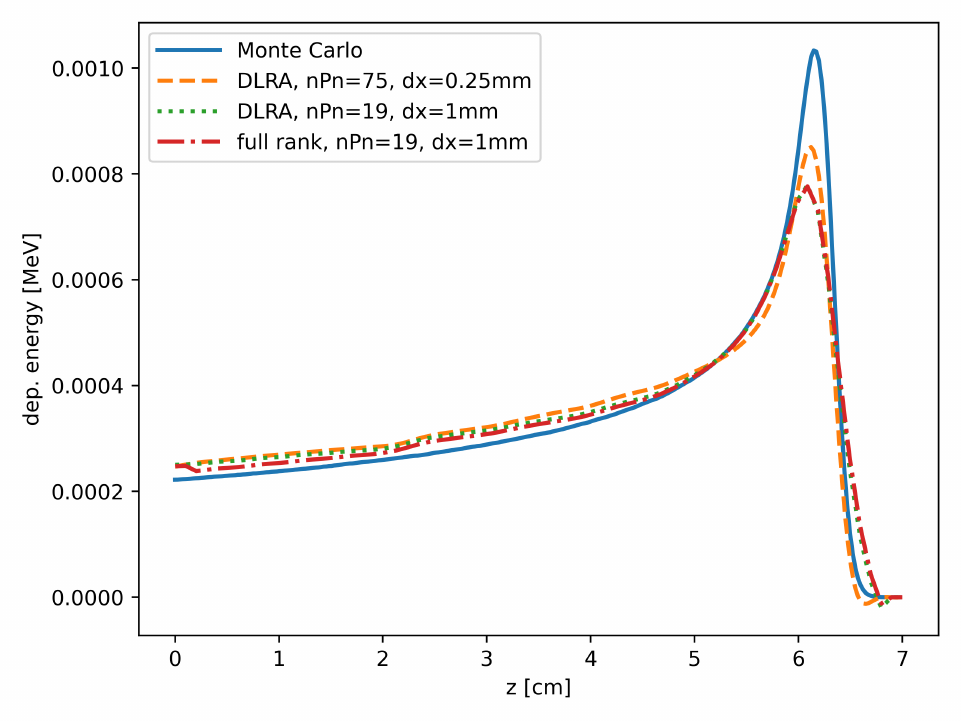}
        \caption{Depth.}
    \end{subfigure}
   \begin{subfigure}{0.31\linewidth}
    \centering
        \includegraphics[width=\linewidth]{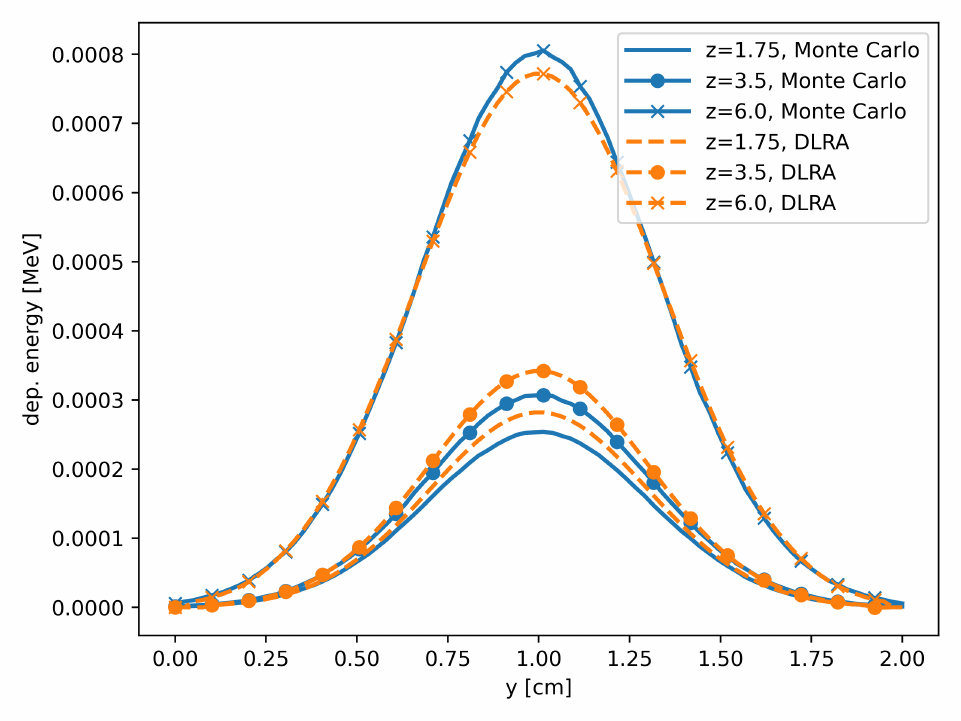}
        \caption{Lateral, MC vs DLRA.}
    \end{subfigure}
 \begin{subfigure}{0.31\linewidth}
  \centering
        \includegraphics[width=\linewidth]{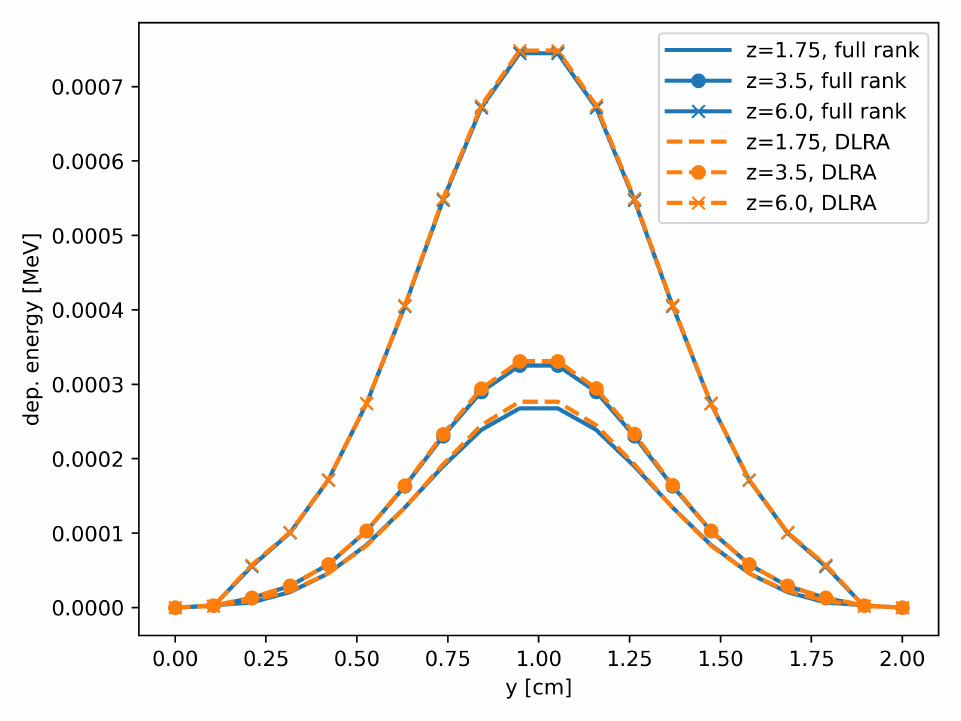}
        \caption{Lateral, DLRA vs full rank.}
    \end{subfigure}
    \caption{Surface plots (a)-(d) and profiles (e)-(g) of deposited energy in the homogeneous test case with a single beam for the Monte Carlo reference  compared to DLRA at high and low spatial and angular resolution and the full-rank deterministic method based on the Fokker-Planck approximation. Lines in (a) indicate positions of profiles plotted in (e)-(g).}
    \label{fig:surfPlotsHomFP}
\end{figure}
   \begin{figure}[h!]
     \centering
     \begin{subfigure}[b]{0.425\textwidth}
         \centering
         \includegraphics[width=\textwidth]{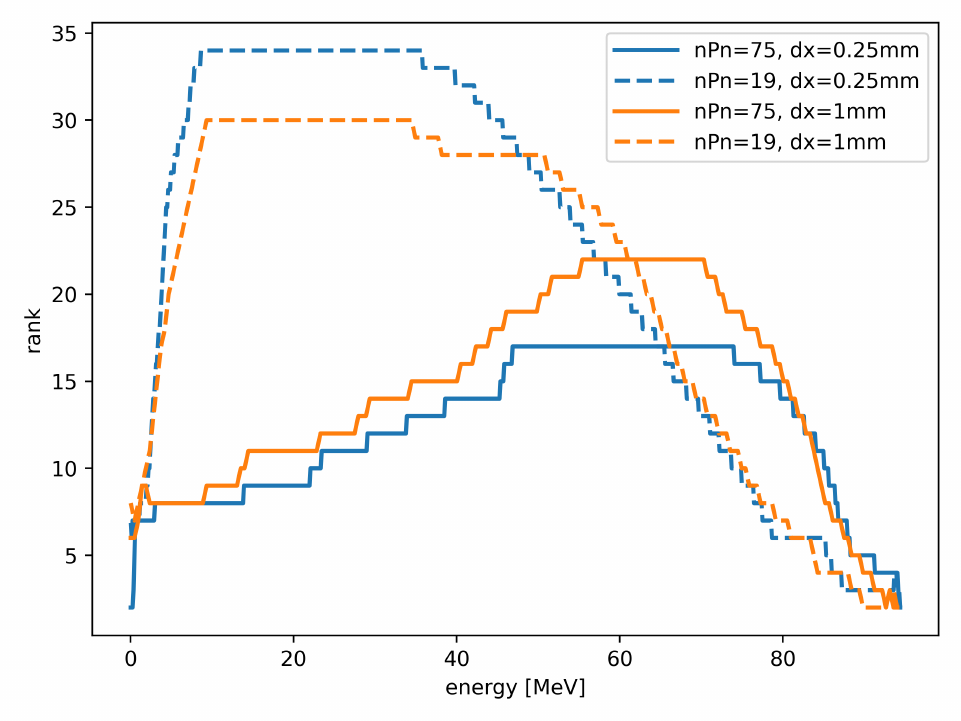}
         \caption{Homogeneous test case.}
     \end{subfigure}
       \begin{subfigure}[b]{0.425\textwidth}
         \centering
         \includegraphics[width=\textwidth]{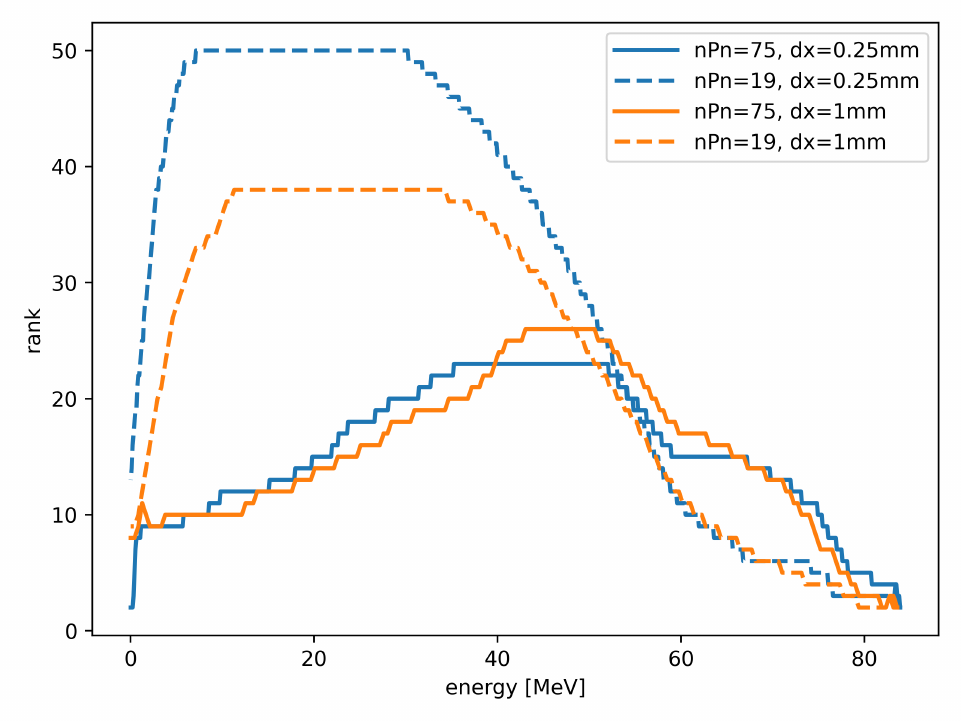}
         \caption{Heterogeneous test case.}
     \end{subfigure}
        \caption{Ranks chosen by rank-adaptive integrator for different spatial and angular resolutions and the Fokker-Planck approximation in (a) the homogeneous and (b) heterogeneous test case with a single beam.}
        \label{fig:RanksFP}
\end{figure} 
 Figure \ref{fig:surfPlotsHomBoltzmann} shows an overview of two-dimensional slices and profiles of the deposited energy using the Boltzmann collision operator at different resolutions compared to the Monte Carlo reference and a full rank deterministic solver. It is apparent, that DLRA can approximate the full-rank solution very well and at a much lower rank according to figure \ref{fig:RanksBoltzmann} (a). It is also clear that the coarse grid resolution is not sufficient to fully capture the characteristics of the solution. When increasing both spatial and angular resolution we can however see that our model is capable of more closely approximating the Monte Carlo reference as well (see figure \ref{fig:surfPlotsHomBoltzmann}). The only significant deviation occurs directly behind the Bragg peak. Here, the solution takes on slightly negative values, which is most likely due to the nature of our angular discretization combined with the steep gradients in this area. Since the full and low-rank solutions match almost exactly, differences to the Monte Carlo reference should however be attributed to the choices of physics models and numerical discretizations.

The Fokker-Planck approximation on the other hand shows smaller differences between the coarse and highly resolved solutions, but tends to underestimate the energy deposition in the region of the Bragg peak and overestimate it slightly in the build-up before the peak as can be seen in figure \ref{fig:surfPlotsHomFP}. This might be attributed to the correction that we use to stabilize the solution \cite{landesman_angular_1989}, which has however proven to be necessary to avoid stronger artifacts near material boundaries during our studies. In figure \ref{fig:surfPlotsHomFP} (e), the solution also exhibits similar negative areas as previously observed for the Boltzmann solver. Figure \ref{fig:RanksFP} further shows that the ranks for the Fokker-Planck approximation are consistently higher than for the Boltzmann equation. Interestingly, in contrast to the Boltzmann solver, the discretizations now affect the rank significantly and finer discretizations are not directly correlated with higher ranks. Here, the ranks are consistently higher for high energies when using a finer angular discretization, but are larger in the case of coarse angular discretizations for small energies. Due to the higher ranks and worse results compared to the Boltzmann solver for high resolutions, using the Fokker-Planck approximation (within the current framework) would only be beneficial for a cheap estimate at a low resolution for which the Boltzmann solver has not converged sufficiently.

\subsection{Heterogeneous Case}
\begin{figure}
    \centering
    \begin{subfigure}{0.225\linewidth}
     \centering
        \includegraphics[width=\linewidth]{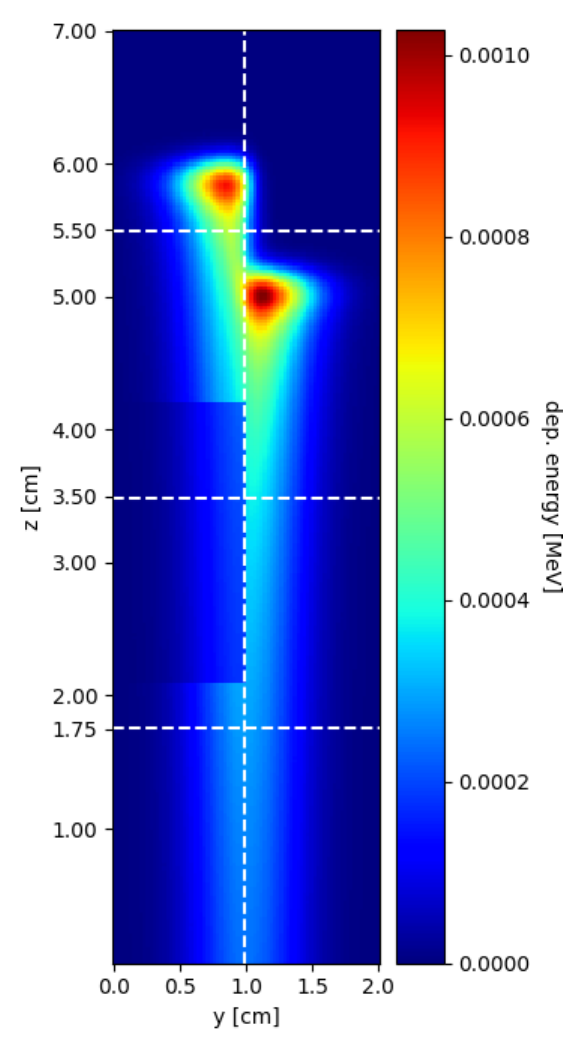}
        \caption{Monte Carlo, \\1e8 histories}
    \end{subfigure}
   \begin{subfigure}{0.225\linewidth}
    \centering
        \includegraphics[width=\linewidth]{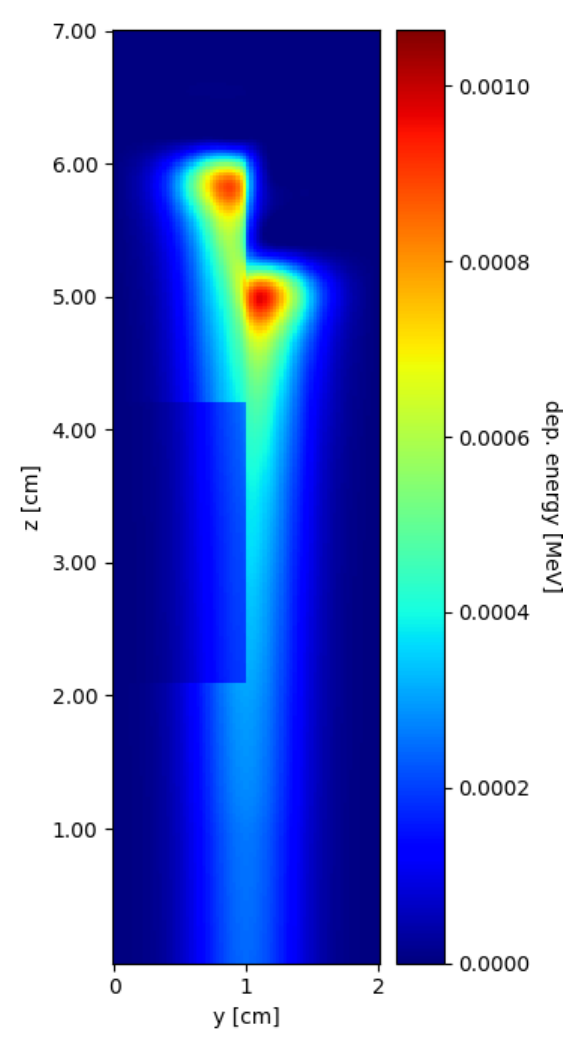}
        \caption{DLRA, P$_{109}$, \\ $\Delta_x=0.25$ mm}
    \end{subfigure}
    \begin{subfigure}{0.225\linewidth}
     \centering
        \includegraphics[width=\linewidth]{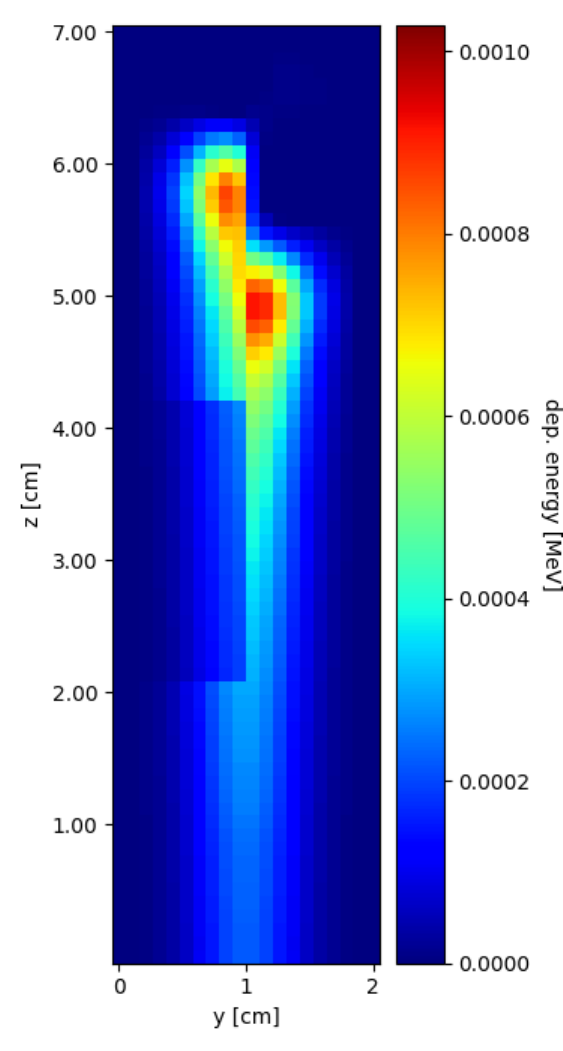}
        \caption{DLRA, P$_{75}$,\\  $\Delta_x=1$ mm}
    \end{subfigure}
    \begin{subfigure}{0.225\linewidth}
     \centering
        \includegraphics[width=\linewidth]{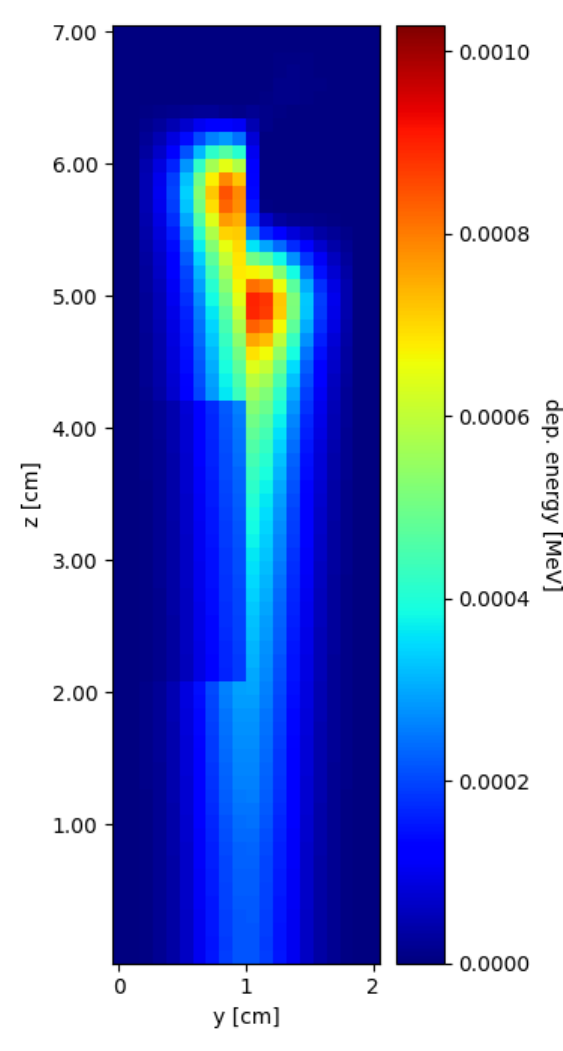}
        \caption{full rank P$_{75}$,\\  $\Delta_x=1$ mm}
    \end{subfigure}\\
        \begin{subfigure}{0.31\linewidth}
         \centering
        \includegraphics[width=\linewidth]{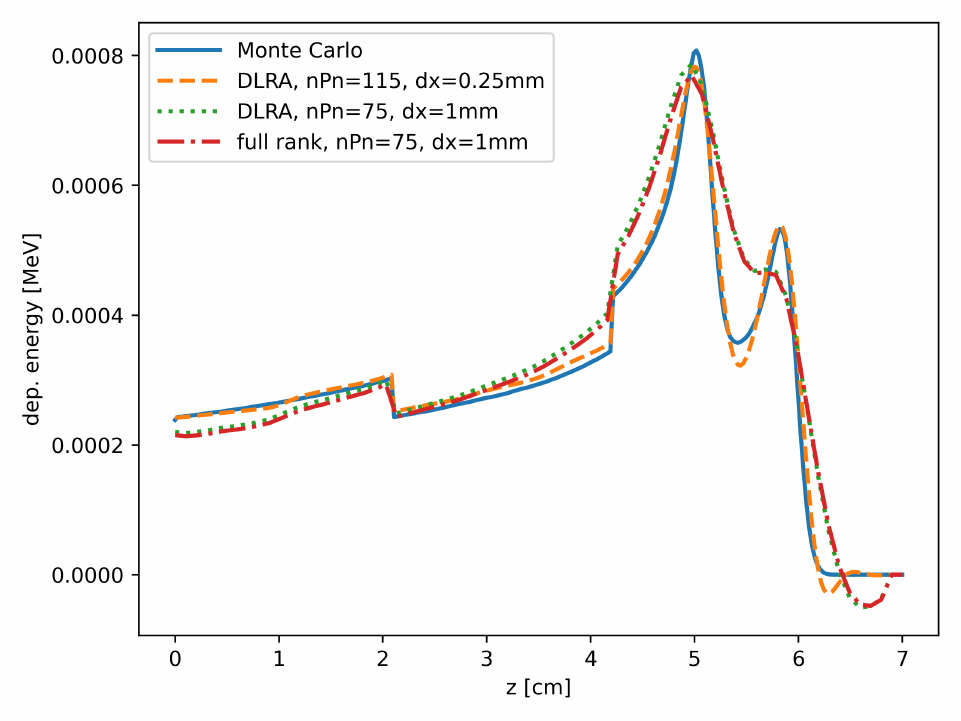}
        \caption{Depth.}
    \end{subfigure}
   \begin{subfigure}{0.31\linewidth}
    \centering
        \includegraphics[width=\linewidth]{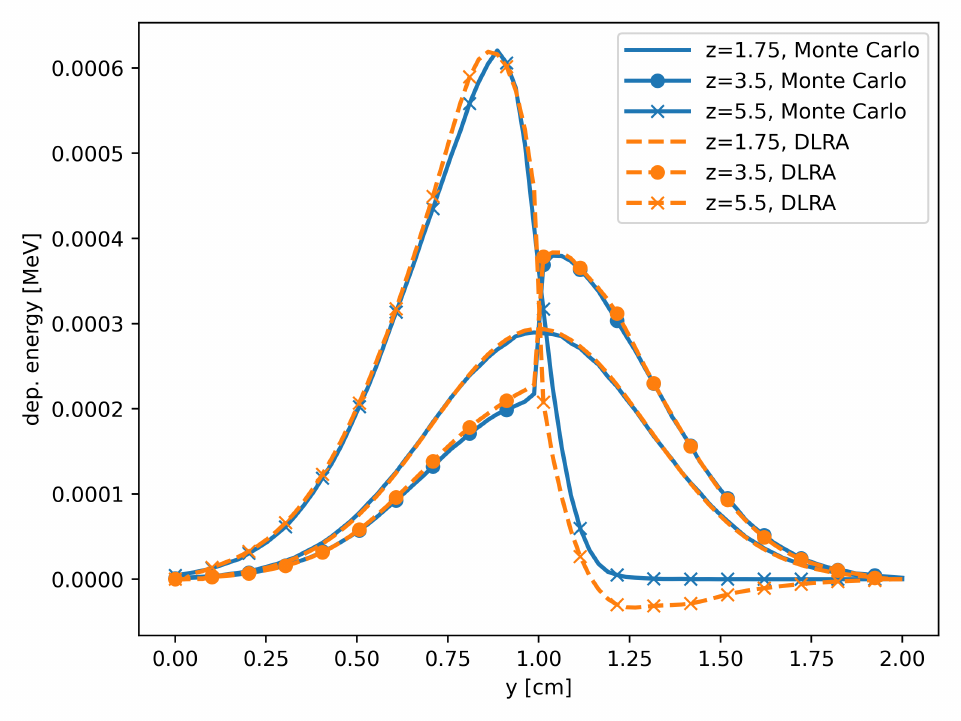}
        \caption{Lateral, MC vs DLRA.}
    \end{subfigure}
 \begin{subfigure}{0.31\linewidth}
  \centering
        \includegraphics[width=\linewidth]{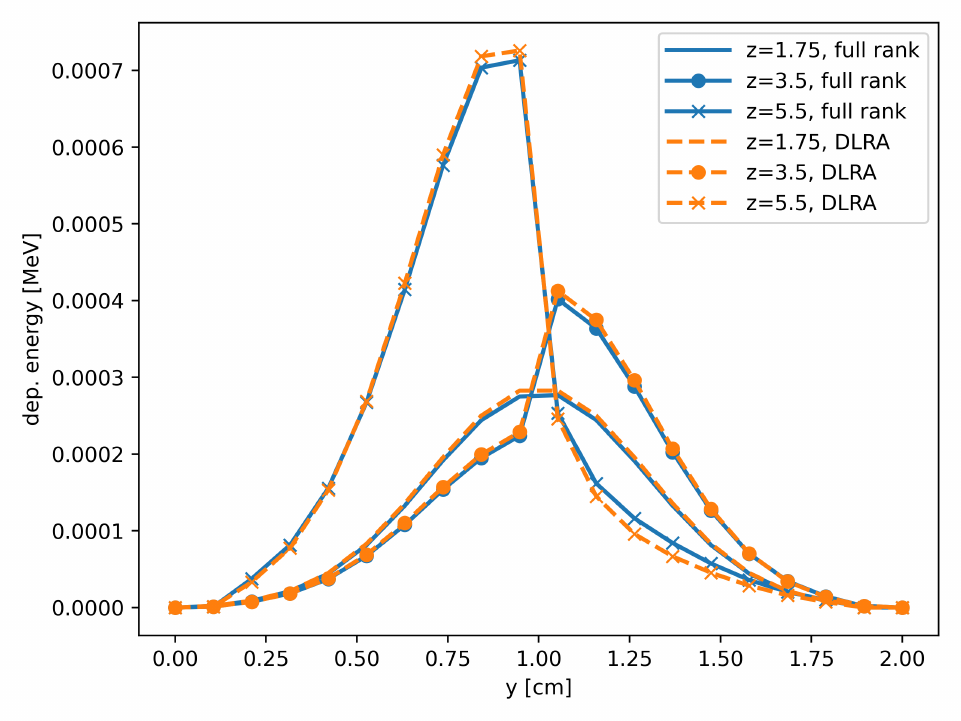}
        \caption{Lateral, DLRA vs full rank.}
    \end{subfigure}
    \caption{Surface plots (a)-(d) and profiles (e)-(g) of deposited energy in the heterogeneous test case with a single beam for the Monte Carlo reference  compared to DLRA at high and low spatial and angular resolution and the full-rank deterministic method based on the Boltzmann equation. Lines in (a) indicate positions of profiles plotted in (e)-(g).}
    \label{fig:SurfPlotsHetBoltzmann}
\end{figure}
Next, we consider a heterogeneous case mimicking a beam encountering strong heterogeneities along its path as might be the case for example in the head and neck region. For this, a lower density box with -400 HU, roughly corresponding to lung tissue, is placed inside the domain (see figure \ref{fig:testCasesB}). A proton beam with 80 MeV is directed at the box such that the center of the beam hits the interface between the materials. Again, the Hounsfield units are converted to material compositions and densities using the method in \cite{schneider_correlation_2000}.
\begin{figure}[h!]
    \centering
        \begin{subfigure}{0.225\linewidth}
         \centering
        \includegraphics[width=\linewidth]{Fig7a_8a.pdf}
        \caption{Monte Carlo, \\1e8 histories}
    \end{subfigure}
   \begin{subfigure}{0.225\linewidth}
    \centering
        \includegraphics[width=\linewidth]{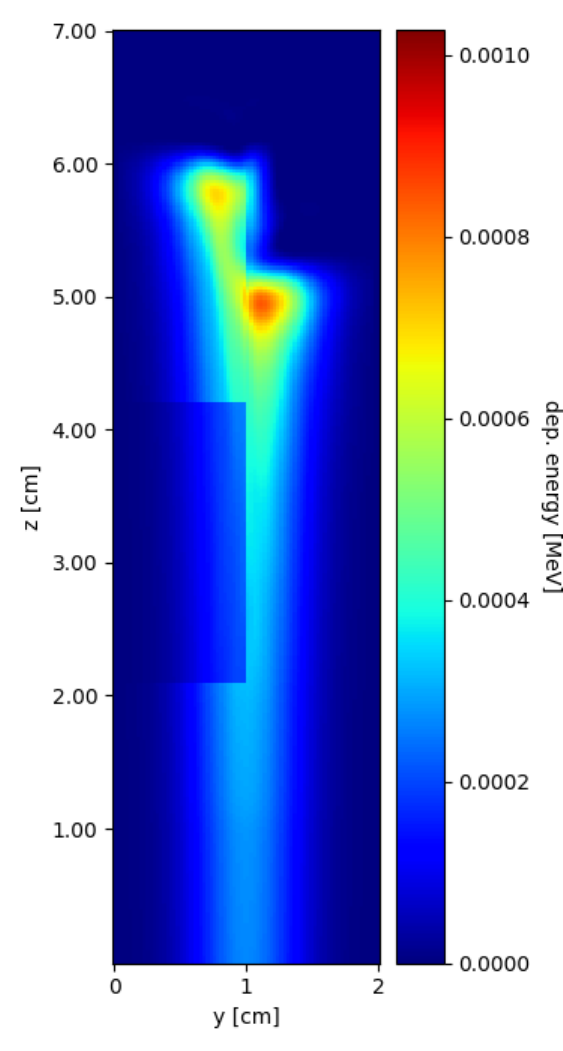}
        \caption{DLRA, P$_{75}$,\\  $\Delta_x=0.25$ mm}
    \end{subfigure}
    \begin{subfigure}{0.225\linewidth}
     \centering
        \includegraphics[width=\linewidth]{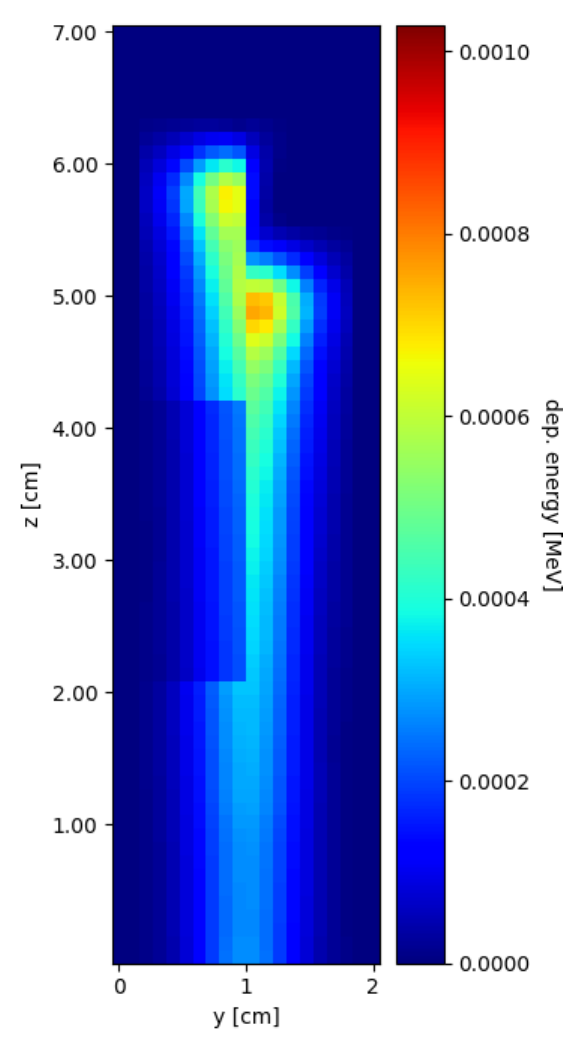}
        \caption{DLRA, P$_{19}$,\\  $\Delta_x=1$ mm}
    \end{subfigure}
    \begin{subfigure}{0.225\linewidth}
     \centering
        \includegraphics[width=\linewidth]{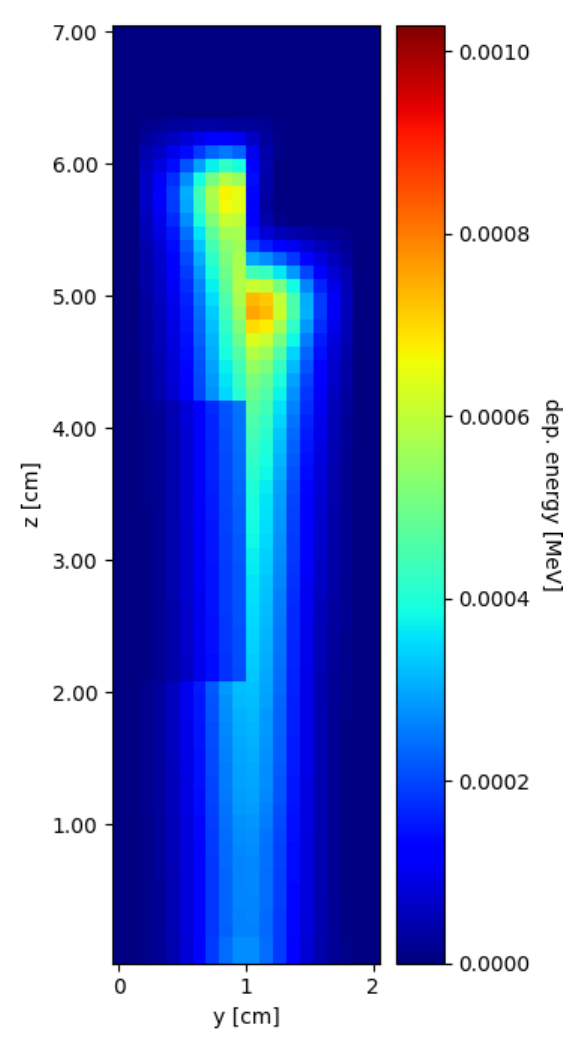}
        \caption{full rank P$_{19}$,\\  $\Delta_x=1$ mm}
    \end{subfigure}\\
            \begin{subfigure}{0.31\linewidth}
             \centering
        \includegraphics[width=\linewidth]{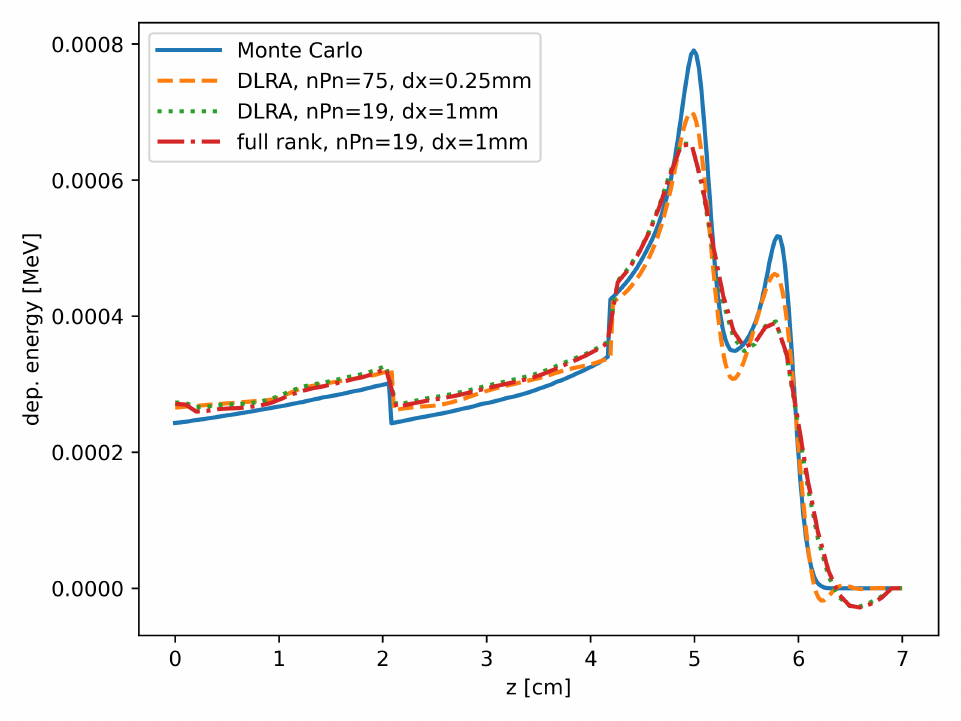}
        \caption{Depth.}
    \end{subfigure}
   \begin{subfigure}{0.31\linewidth}
    \centering
        \includegraphics[width=\linewidth]{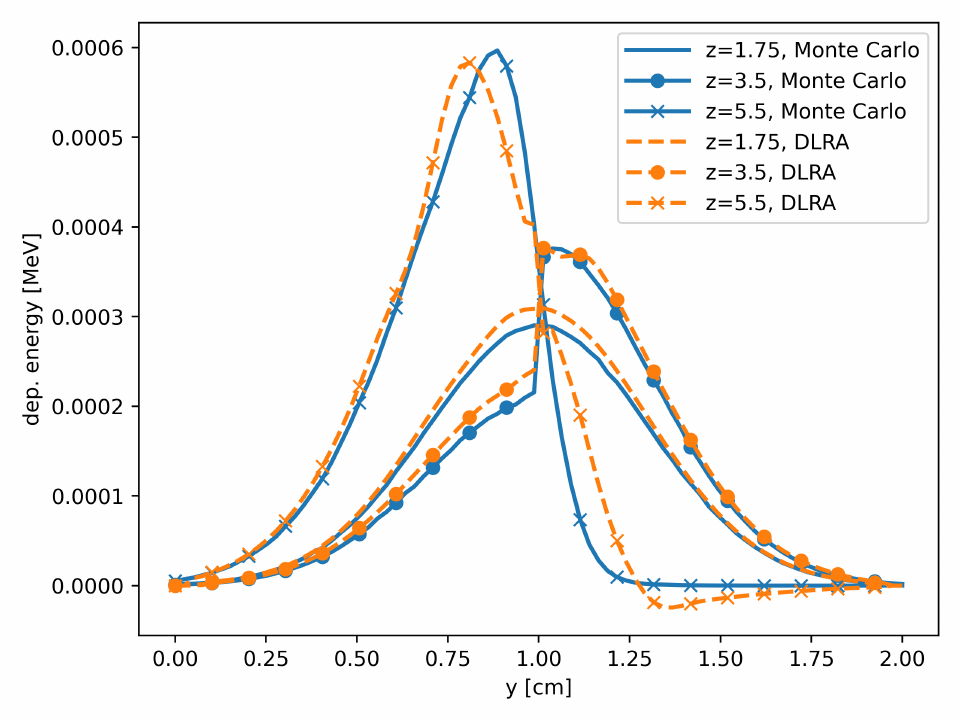}
        \caption{Lateral, MC vs DLRA.}
    \end{subfigure}
 \begin{subfigure}{0.31\linewidth}
  \centering
        \includegraphics[width=\linewidth]{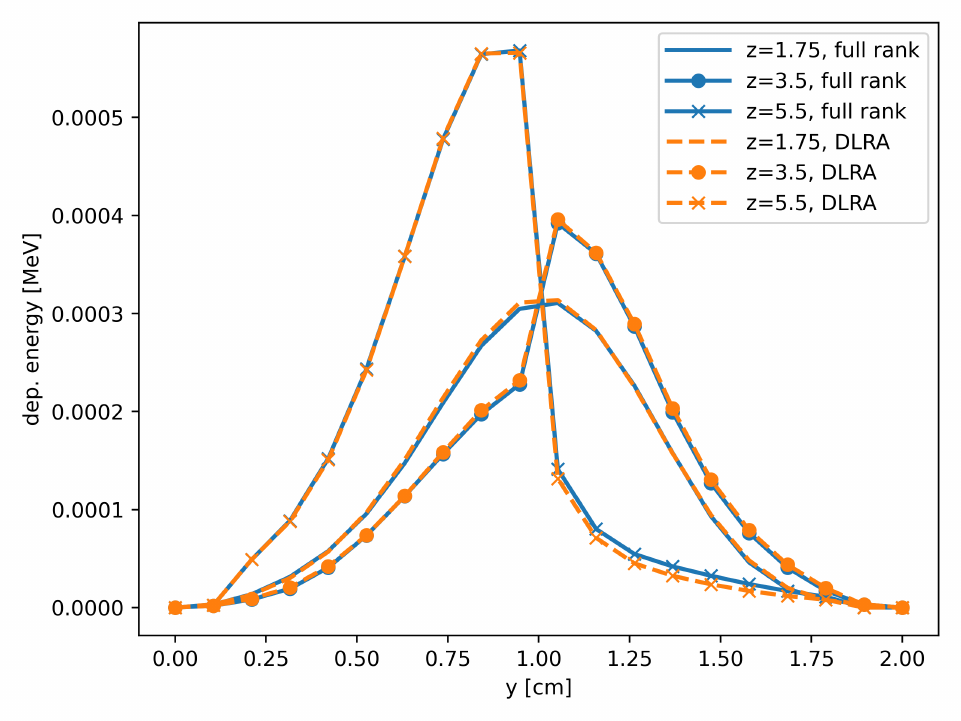}
        \caption{Lateral, DLRA vs full rank.}
    \end{subfigure}
    \caption{Surface plots (a)-(d) and profiles (e)-(g) of deposited energy in the heterogeneous test case with a single beam for the Monte Carlo reference  compared to DLRA at high and low spatial and angular resolution and the full-rank deterministic method based on the Fokker-Planck approximation. Lines in (a) indicate positions of profiles plotted in (e)-(g).}
    \label{fig:surfPlotsHetFP}
\end{figure}
As in the homogeneous case, the full-rank solution is well-approximated by the dynamical low-rank approximation for both the Boltzmann and the Fokker-Planck solvers. This can be seen both in the 2D surface plots in figure \ref{fig:SurfPlotsHetBoltzmann} (c) and (d) (Boltzmann) and figure \ref{fig:surfPlotsHetFP} (c) and (d) (Fokker-Planck) and the one-dimensional solution profiles in figure \ref{fig:SurfPlotsHetBoltzmann}(e) and (g). The Boltzmann solver also approximates the Monte Carlo reference reasonably well, although more deviations can be seen than in the homogeneous test case. For example, the surface plots in figure \ref{fig:SurfPlotsHetBoltzmann}(a) and (b) show that the highest dose values in the Bragg peaks are slightly underestimated and depth profiles in figure \ref{fig:SurfPlotsHetBoltzmann} (e) reveal differences in the local minimum behind the first peak. The ranks required for this more complex test case are expectedly higher than in the homogeneous case (see figure \ref{fig:RanksBoltzmann} (b)) and increase for increasingly fine spatial and angular resolutions. All chosen ranks are however still significantly smaller than the number of spatial cells or angular modes and lead to a large reduction of the system size (see also section \ref{sec:complexity}). 

When comparing the Fokker-Planck approximation to the Monte Carlo reference in figures \ref{fig:surfPlotsHetFP} (a) and (b) we observe similar results as in the homogeneous case as well. Again, the solution tends to underestimate the highest energy depositions and there are small negative areas behind the Bragg peak. In this case, the sharp material boundary however tended to cause artifacts, which can still be slightly seen in figure \ref{fig:surfPlotsHetFP} (b). While a stronger correction of the angular moments according to \cite{landesman_angular_1989,morel_fokker-planck_1981} reduces the artifacts, it also has a smoothening effect and dampens the height of the Bragg peak, creating a tradeoff between stability and preserving solution characteristics. Figure \ref{fig:RanksFP} (b) further shows that ranks for different discretizations have a similar structure as in the homogeneous case, but are consistently higher, as already observed for the Boltzmann solver.

\subsection{Homogeneous with two beams}
\begin{figure}[h!]
    \centering
     \rotatebox{90}{\hspace{1cm}\textbf{Boltzmann}}
    \begin{subfigure}{0.235\linewidth}
        \includegraphics[width=\linewidth]{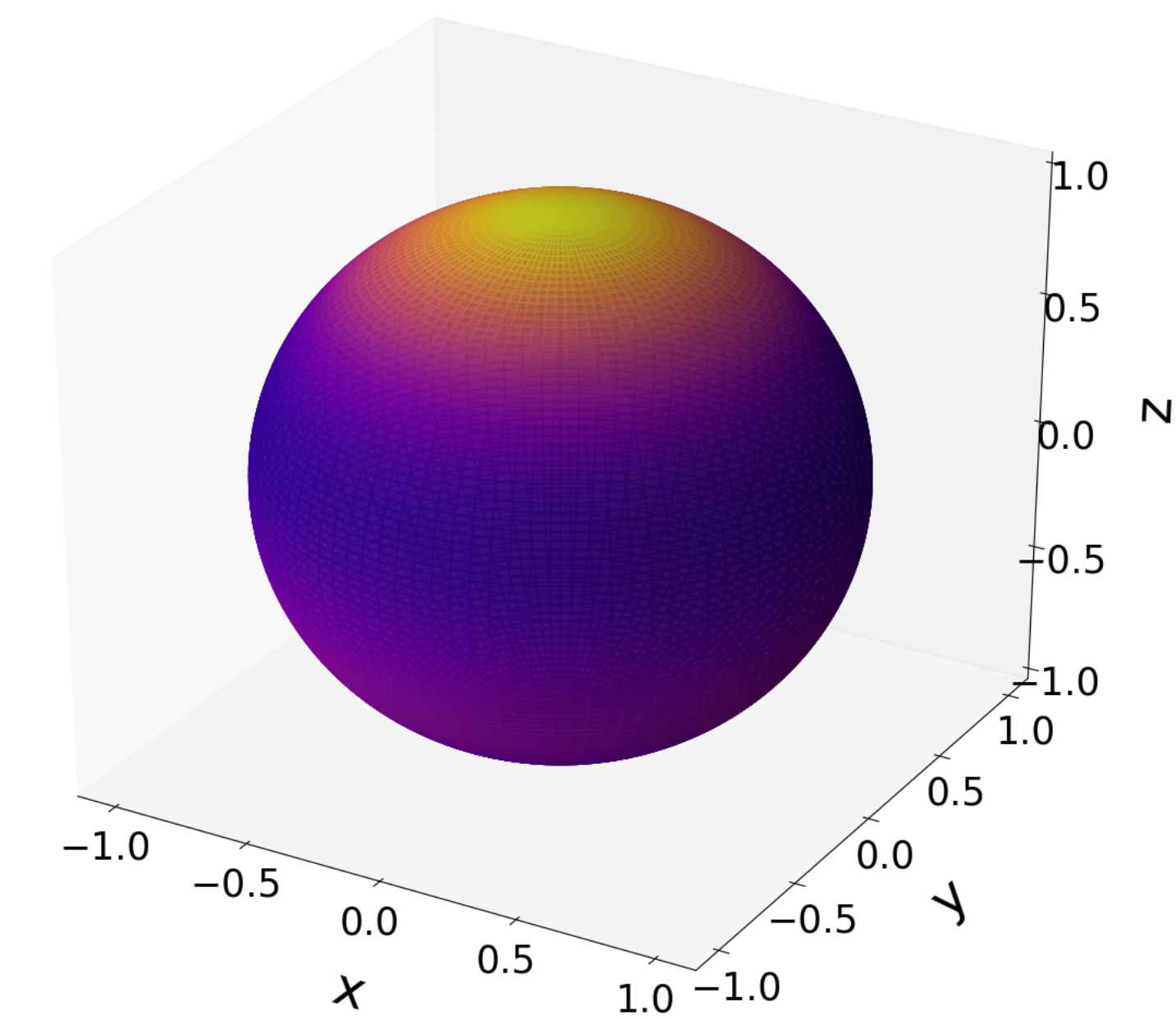}
        \caption{Mode 1}
    \end{subfigure}
    \begin{subfigure}{0.235\linewidth}
        \includegraphics[width=\linewidth]{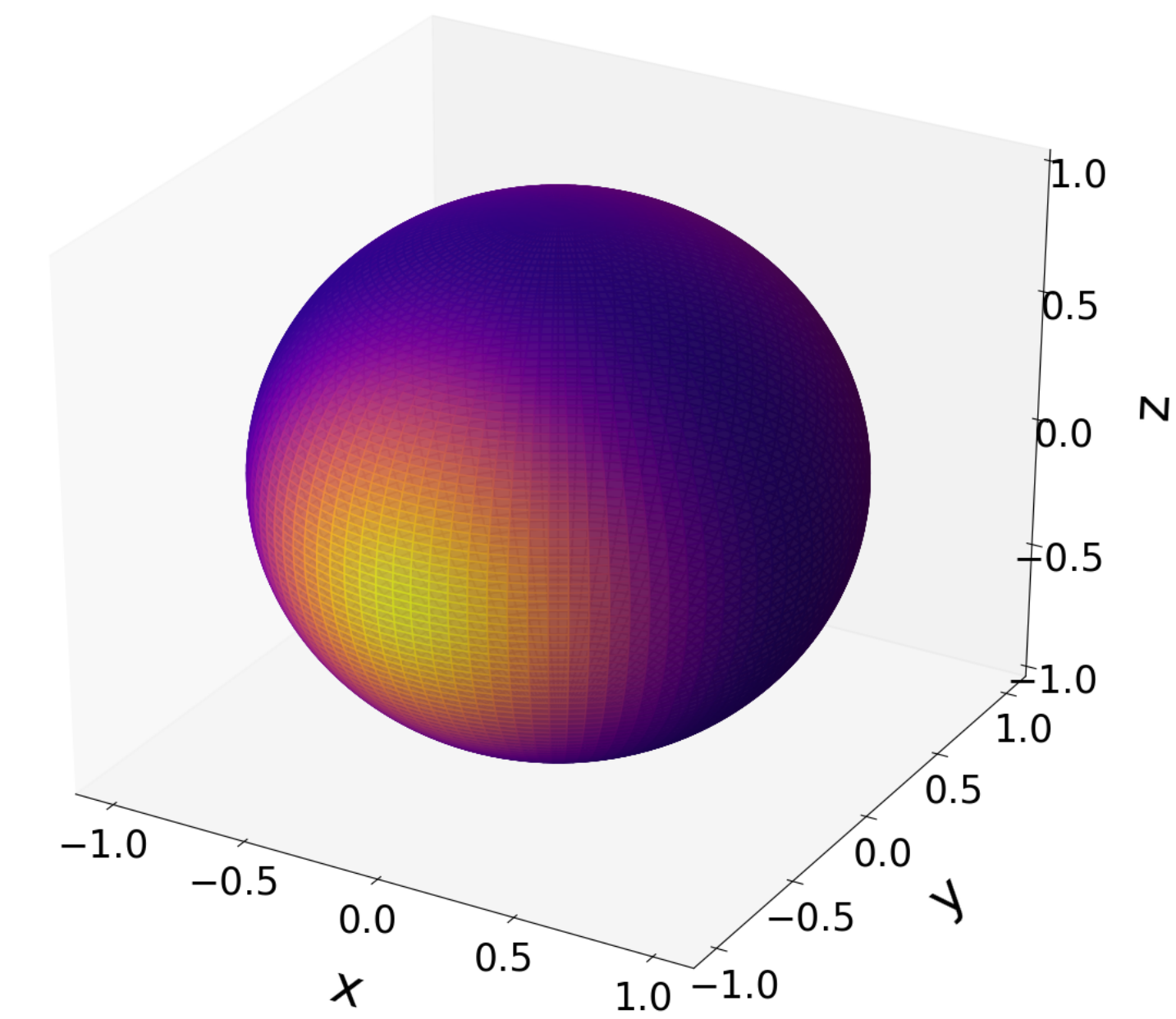}
        \caption{Mode 2}
    \end{subfigure}
     \begin{subfigure}{0.235\linewidth}
        \includegraphics[width=\linewidth]{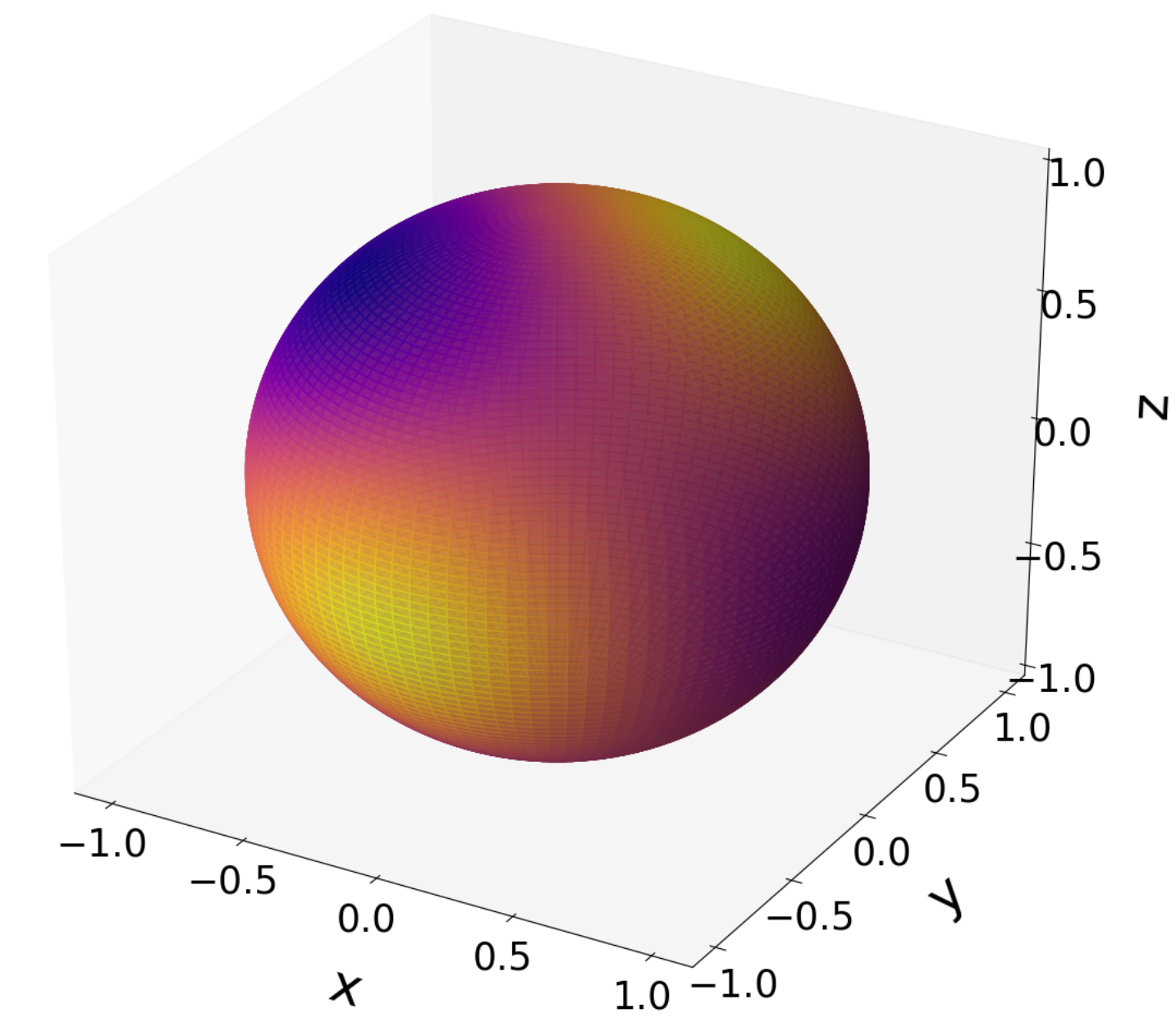}
        \caption{Mode 3}
    \end{subfigure}
    \begin{subfigure}{0.235\linewidth}
        \includegraphics[width=\linewidth]{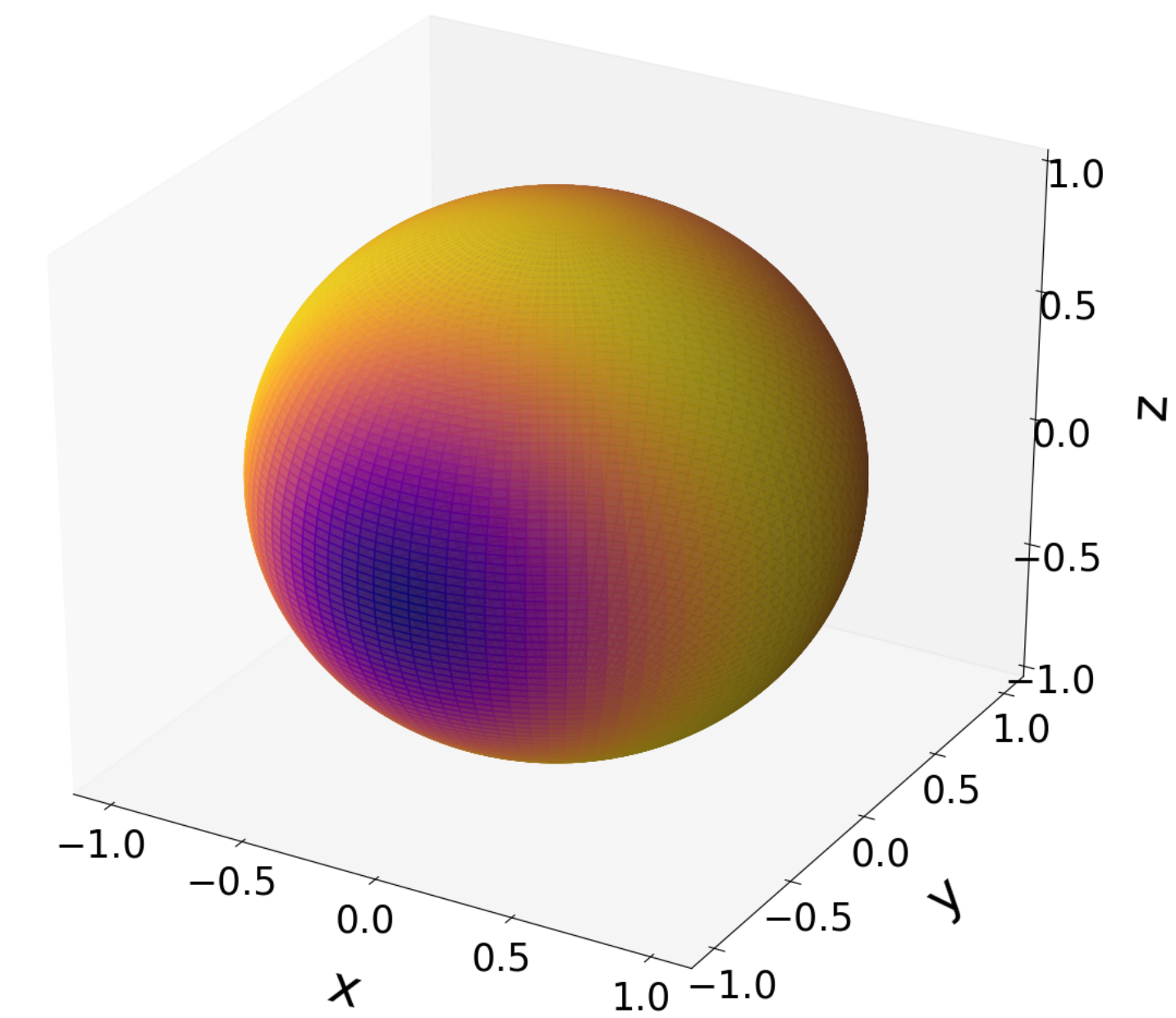}
        \caption{Mode 4}
    \end{subfigure}\\
     \rotatebox{90}{\hspace{1cm}\textbf{Fokker-Planck}}
         \begin{subfigure}{0.235\linewidth}
        \includegraphics[width=\linewidth]{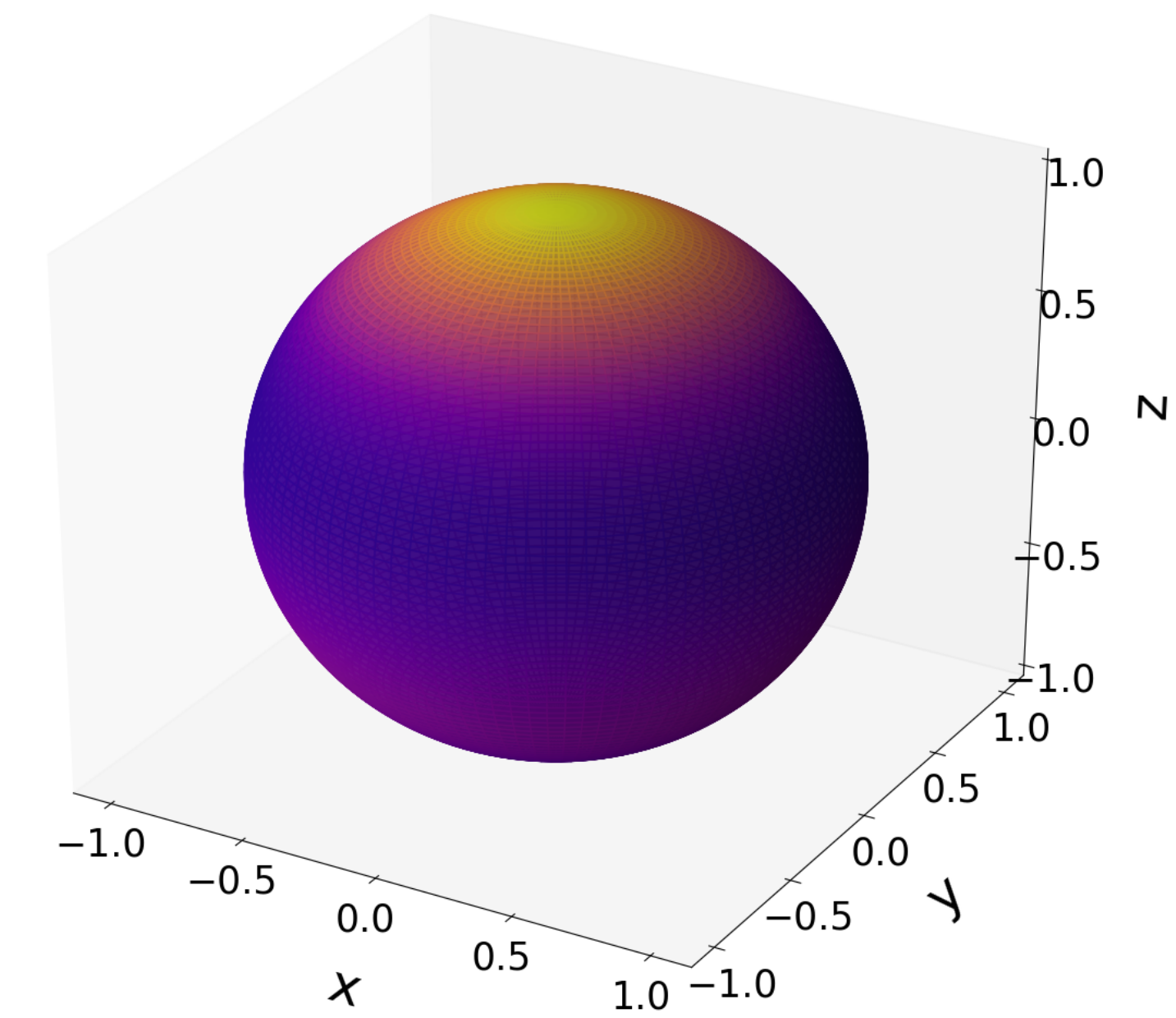}
        \caption{Mode 1}
    \end{subfigure}
    \begin{subfigure}{0.235\linewidth}
        \includegraphics[width=\linewidth]{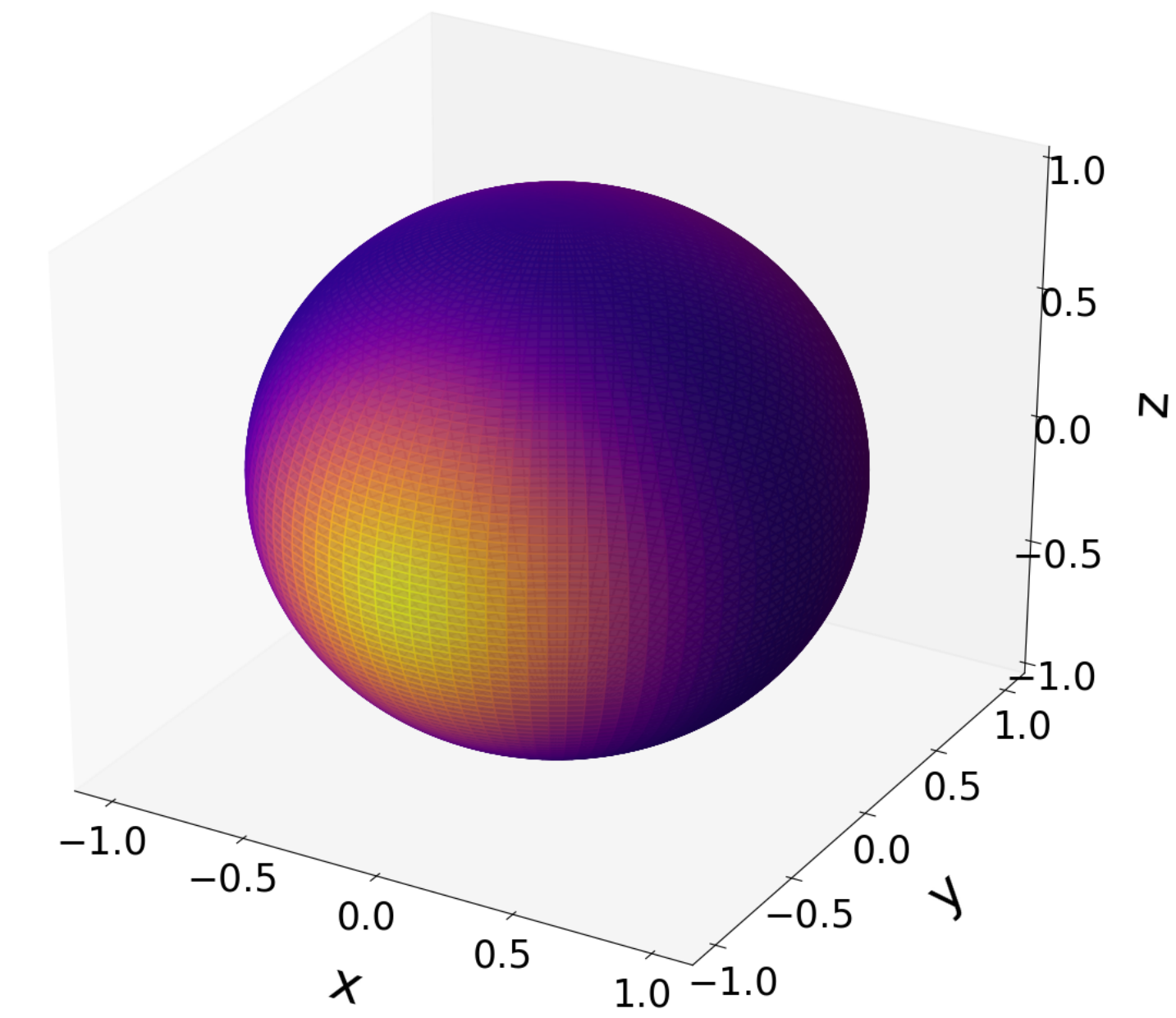}
        \caption{Mode 2}
    \end{subfigure}
        \begin{subfigure}{0.235\linewidth}
        \includegraphics[width=\linewidth]{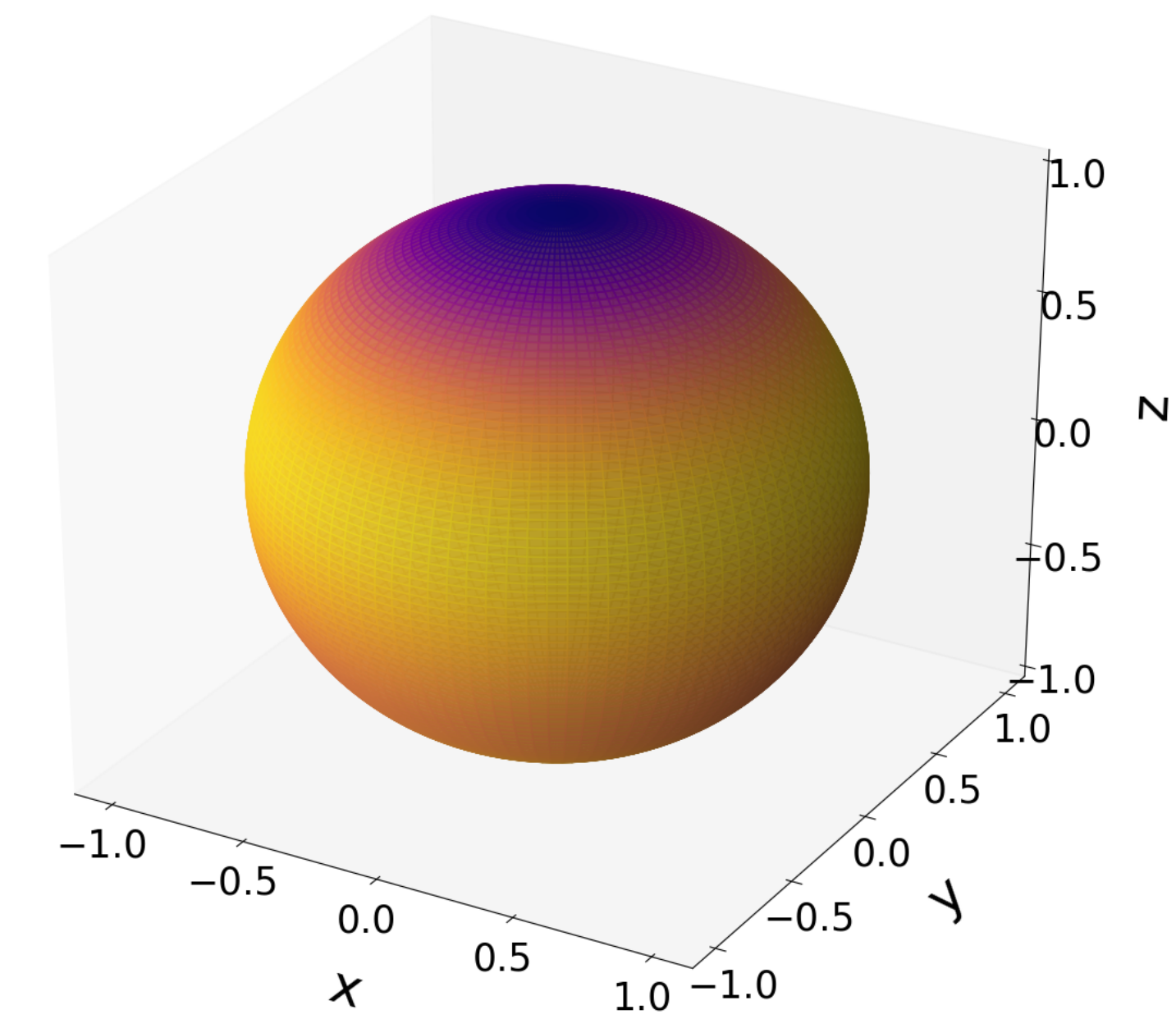}
        \caption{Mode 3}
    \end{subfigure}
    \begin{subfigure}{0.235\linewidth}
        \includegraphics[width=\linewidth]{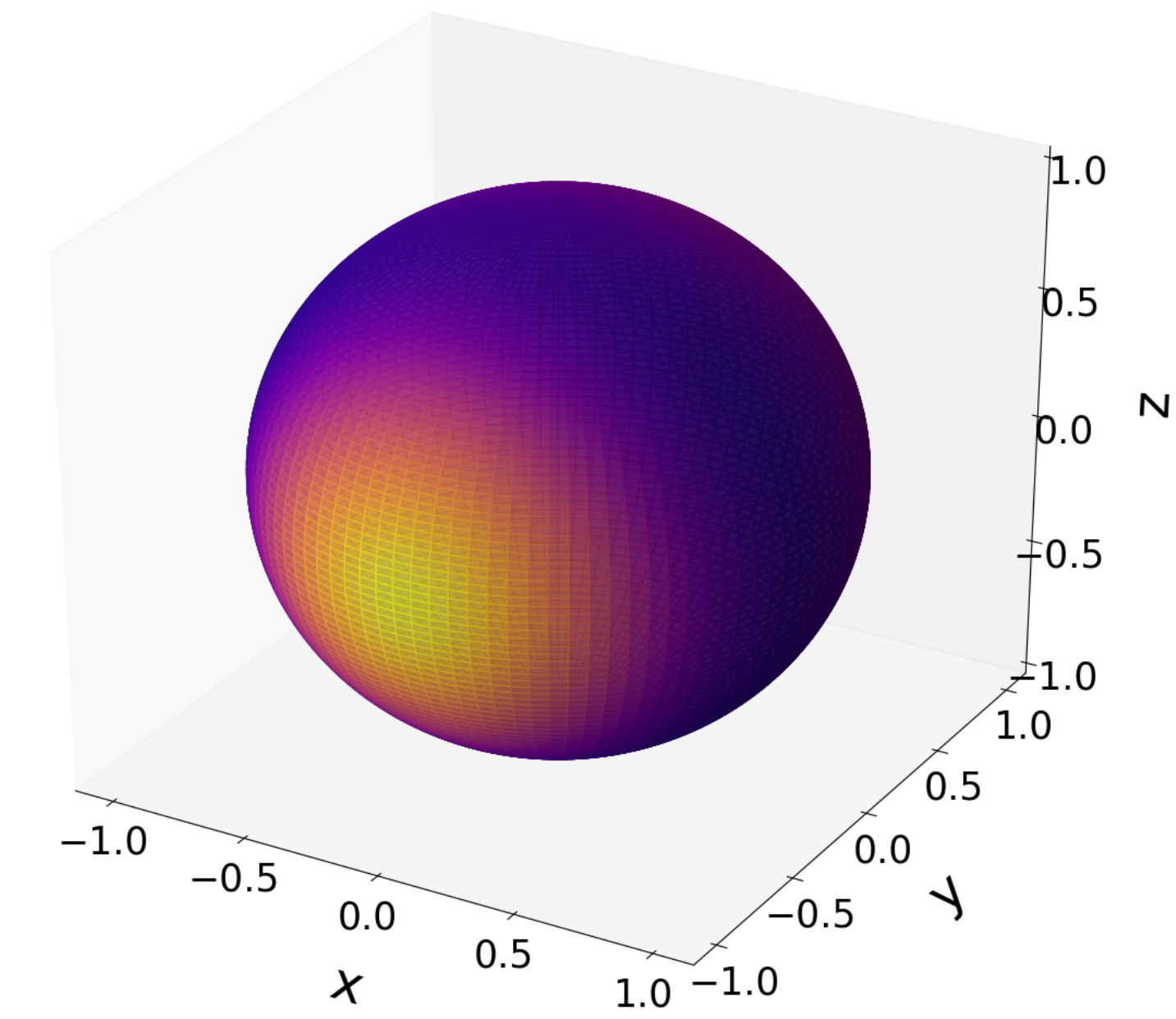}
        \caption{Mode 4}
    \end{subfigure}
    \caption{First four dominant angular modes after one fifth of the time/energy steps.}
    \label{fig:angBasisTwoBeams}
\end{figure}

Finally, we consider a case with two 50 MeV beams entering the domain at different angles (see figure \ref{fig:testCasesC}). A treatment plan in proton therapy typically consists of a larger number of layered beamlets with different initial energies, positions and angles. Thus, this is a good test case to gauge the impact of a more complex beam source on rank and therefore computational costs and memory. 
\begin{figure}[h!]
    \centering
    \begin{subfigure}{0.325\linewidth}
        \includegraphics[width=\linewidth]{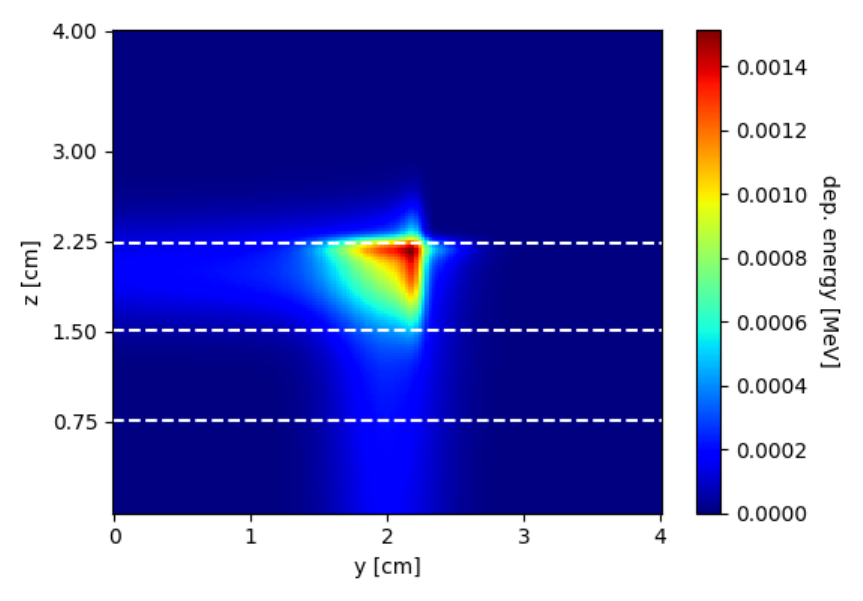}
        \caption{Monte Carlo}
    \end{subfigure}
    \begin{subfigure}{0.325\linewidth}
        \includegraphics[width=\linewidth]{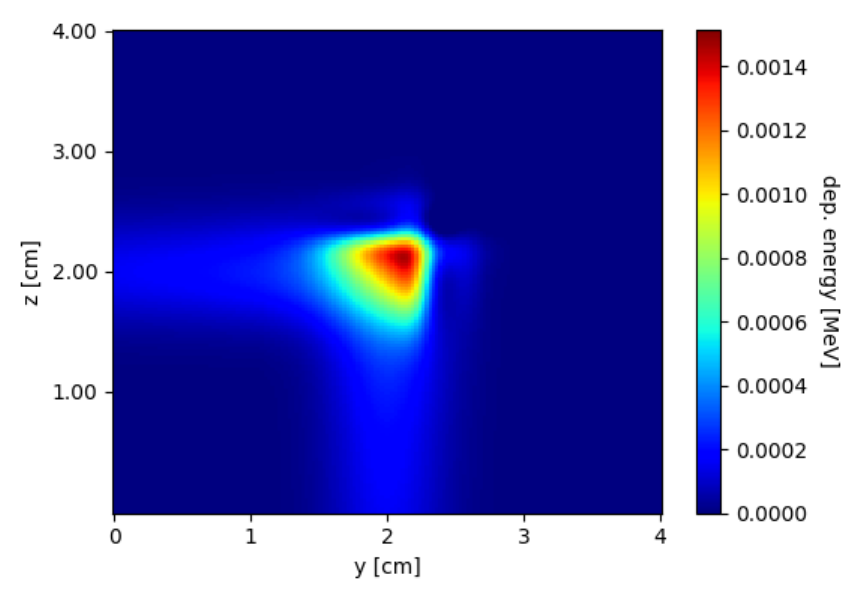}
        \caption{DLRA, Boltzmann}
    \end{subfigure}
       \begin{subfigure}{0.325\linewidth}
        \includegraphics[width=\linewidth]{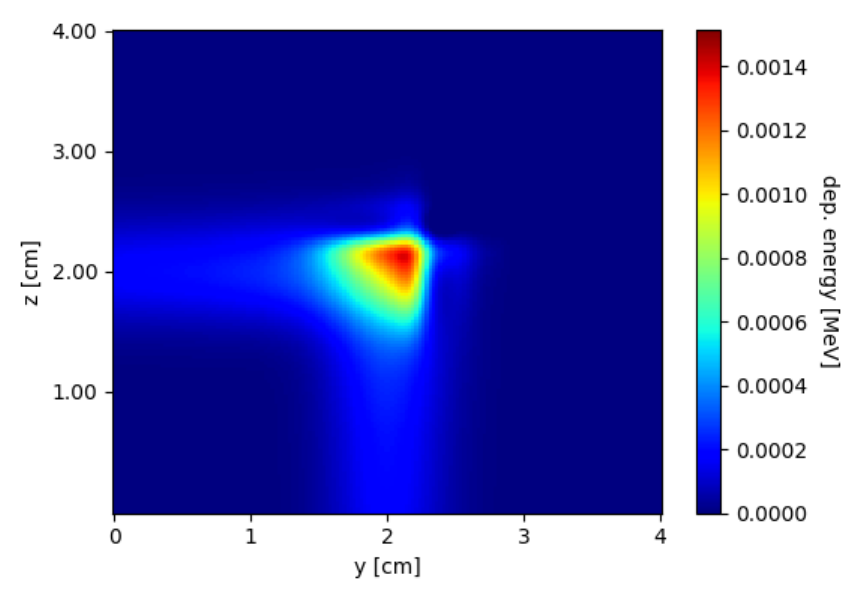}
        \caption{DLRA, Fokker-Planck}
    \end{subfigure}\\
     \begin{subfigure}{0.425\linewidth}
    \centering
        \includegraphics[width=\linewidth]{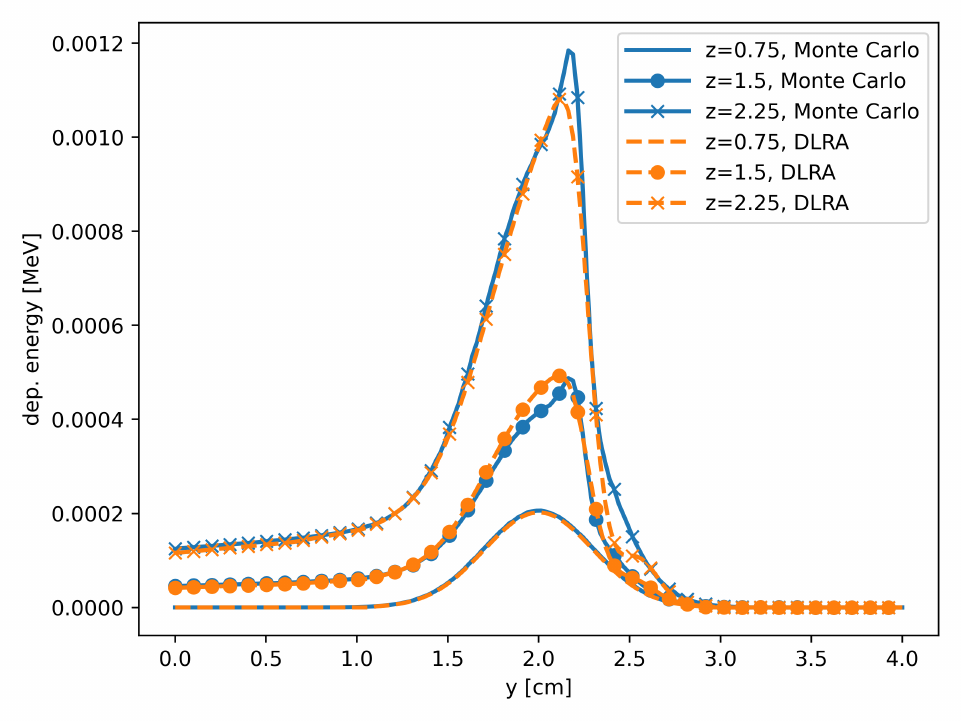}
        \caption{Profiles, Boltzmann.}
    \end{subfigure}
       \begin{subfigure}{0.425\linewidth}
    \centering
        \includegraphics[width=\linewidth]{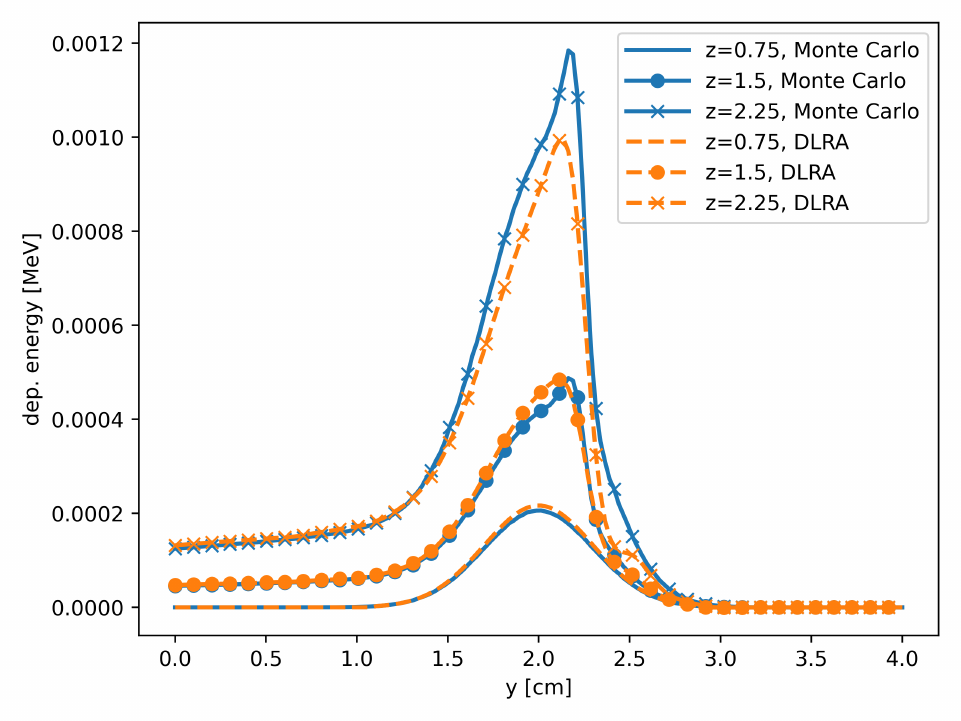}
        \caption{Profiles, Fokker-Planck.}
    \end{subfigure}
       \caption{Surface plots (a)-(c) and profiles (d)-(e) of deposited energy for two perpendicular beams in a homogeneous material for the Monte Carlo reference compared to DLRA based on the Boltzmann equation and the Fokker-Planck approximation. Lines in (a) indicate positions of profiles plotted in (d)-(e).}
        \label{fig:resTwoBeams}
    \end{figure}
    \begin{figure}[h!]
        \centering
     \begin{subfigure}{0.425\linewidth}
      \centering
        \includegraphics[width=\linewidth]{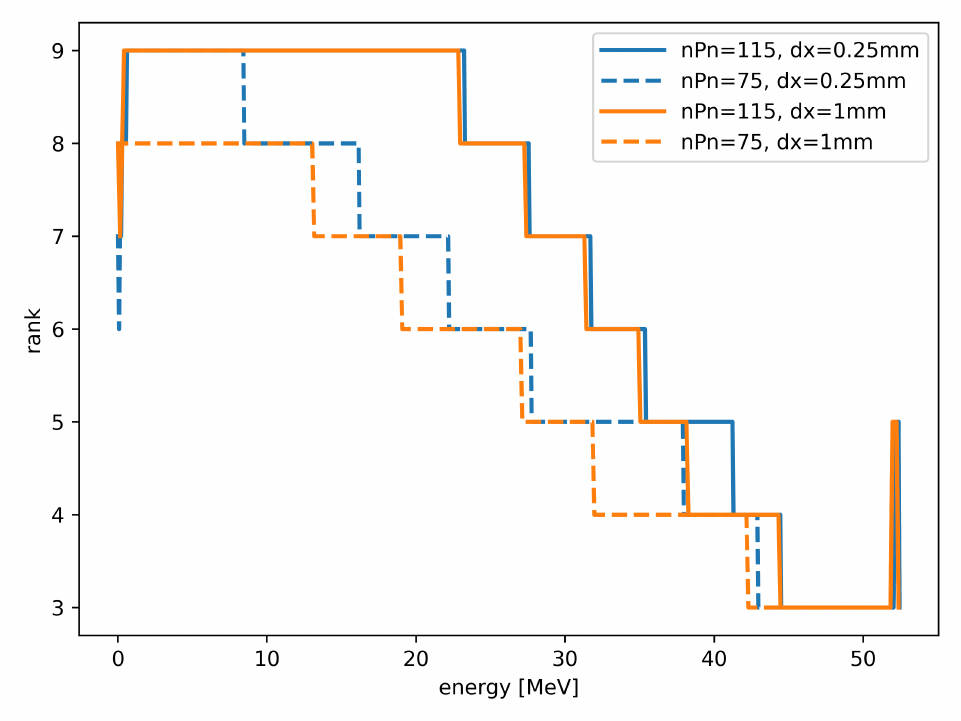}
        \caption{Boltzmann.}
    \end{subfigure}
     \begin{subfigure}{0.425\linewidth}
      \centering
        \includegraphics[width=\linewidth]{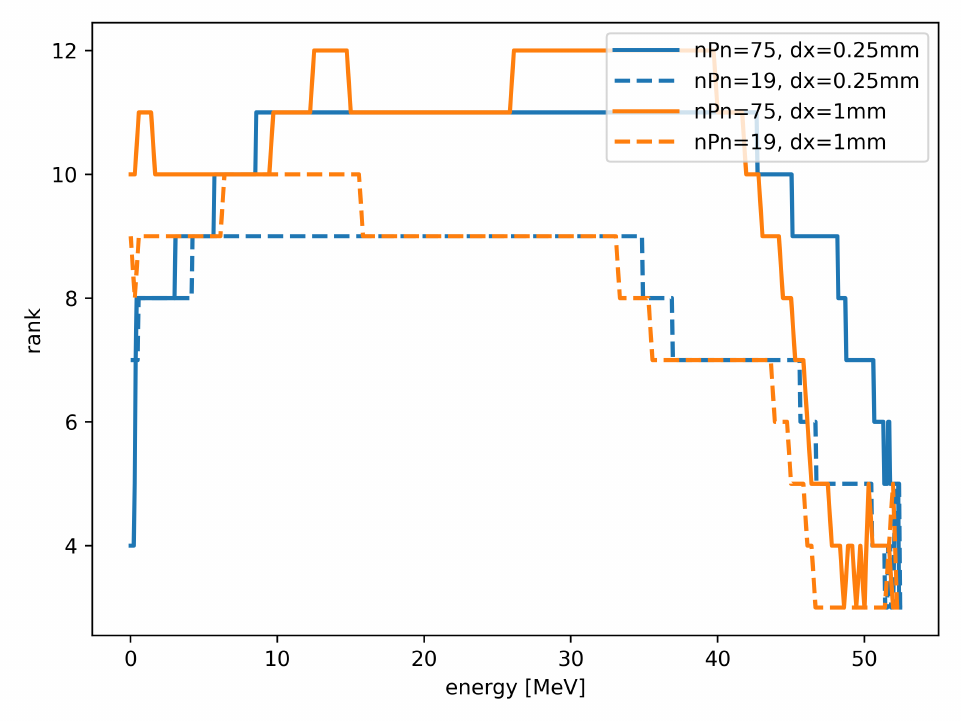}
        \caption{Fokker-Planck.}
    \end{subfigure}
    \caption{Ranks chosen by rank-adaptive integrator using (a) Boltzmann and (b) the Fokker-Planck approximation in the homogeneous test case with two beams.}
    \label{fig:ranksTwoBeams}
\end{figure}
Figures \ref{fig:resTwoBeams} and \ref{fig:ranksTwoBeams} show a comparison of Monte Carlo to DLRA as well as the ranks chosen using the same tolerance as in the previous test cases. The low-rank method manages to capture the basic dose characteristics well, but appears to slightly smoothing in the sense of losing local "hot-spots" in the solution and corresponding sharp dose gradients. For the Boltzmann solver, the chosen ranks are on average roughly 2 points higher than in the homogeneous single beam test case for the finer angular resolution and even less when comparing the coarser resolutions (compare figures \ref{fig:RanksBoltzmann} (a) and \ref{fig:ranksTwoBeams} (a)). Interestingly, comparing figures \ref{fig:RanksFP} (b) and \ref{fig:ranksTwoBeams} (b), we see that the ranks are even lower than those of the single beam test case for the Fokker-Planck approximation. This difference could also be due to the lower beam energy in this test case. While numerical experiments have shown that for a tolerance of $\vartheta=0.01$ the initial beam energy does not affect the rank for the Boltzmann solver, we have seen that the ranks for the Fokker-Planck solver decreased for decreasing beam energies. Further, for both approaches the ranks are mostly affected by the angular discretization, which is somewhat intuitive considering we effectively doubled the amount of relevant directions.  Figure \ref{fig:angBasisTwoBeams} also shows that for both Boltzmann and Fokker-Planck, the first two modes of the angular basis indeed represent the directions of the two beams, thus capturing the forward movement of particles with some angular spread due to scattering. 

Since in theory the rank quadratically impacts the computational costs (see section \ref{sec:complexity}), further studies are necessary to see when it is more efficient to treat beam sources separately or jointly. Our results however suggest that at the very least a relatively narrow spread in angle, which is typically used for beamlets originating from the same beam direction, can be incorporated without a significant increase in computational costs or memory. This would already reduce the amount of individual computations to only the 2-5 beam angles, compared to potentially thousands of pencil beams \cite{cao_reflections_2022}.

\section{Memory, computational complexity and implementation}
\label{sec:complexity}
The dynamical low-rank approximation reduces the computational complexity from $\mathcal{O}(n\cdot m)$ to $\mathcal{O}(r^2\cdot(n+m))$ and $\mathcal{O}(r\cdot(n+m))$ with respect to runtime and memory, respectively \cite{kusch_robust_2023}. Thus, for $r\ll m, n $ we can expect significant cost reductions. In order to conduct our numerical experiments, the described method is implemented in julia with GPU acceleration using CUDA.jl. The code to reproduce our numerical experiments is openly available on \href{https://github.com/piastammer/publication-DLRA-for-Proton-Transport}{GitHub}. Since the full-rank method with a meaningful discretization does not fit on the GPU, we further implement a full-rank CPU version using the same numerical methods and physical models. All experiments are run on a workstation PC\footnote{CPU: Dual Intel Xeon Silver 4214R (24 cores, 48 threads) ~2.9GHz, GPU: 3 NVIDIA A40, 46GB memory, 1.15–1.74GHz (only 1 GPU used here)} which is comparable to what can be expected in a clinical setting.

Table \ref{tab:timememory} shows a comparison of the full-rank and DLRA Boltzmann and Fokker-Planck solvers with respect to their theoretical as well as measured costs\footnote{The reported memory is total allocated memory according to TimerOutputs.jl. To avoid ambiguity we use GiB (Gibibyte) and TiB (Tebibyte) to denote the memory in the binary system, i.e., 1 GiB = $1024^3$ bytes. Times are recorded using the same julia package with DLRA running on the GPU and the full rank method on the CPU with only the inherent linear algebra parallelization from julia.} for the homogeneous test case with coarser spatial and angular resolution. We report the minimum excecution time out of 20 runs, since the execution path is completely deterministic and thus variations are due to system overhead and background processes on the shared workstation. The theoretical values are computed using the formulas for computational complexity and memory from above and the average rank over time/energy steps is used for DLRA, which is 4.77 and 19.99 for Boltzmann and Fokker-Planck, respectively. We see that in case of the Boltzmann solver, DLRA reduces the costs to less than 0.5\% in all cases and consistently achieves an even higher efficiency gain in practice, which demonstrates the combined effects of GPU implementation and complexity reduction. As predicted theoretically, the reductions are larger in terms of memory consumption. 

Since for the Fokker-Planck approximation we both use lower orders of angular discretization and the average rank is higher than for Boltzmann, there are no theoretical complexity improvements in the case considered here. In practice, we can however see almost two orders of magnitude speed-up due to the GPU implementation. Thus, even for slightly higher rank and less high-dimensional problems, DLRA can be worthwhile by reducing the memory. Further, we see in figure \ref{fig:RanksFP} that the ranks for Fokker-Planck do not necessarily increase for finer spatial discretizations, so for more highly resolved problems we would also see larger theoretical gains with respect to runtime.
 \begin{table}[h!]
        \centering
        \caption{Theoretical and measured time and memory costs of DLRA vs. full rank.}
        \label{tab:timememory}
        \begin{tabular}{c c c | c }
        \textit{Boltzmann}    & full rank (CPU) & DLRA (GPU) & \% of full rank \\
            \hline 
          complexity & \num{161728000} & \num{769597.8992} &0.4759\%\\
          system size& \num{161728000} & \num{161226.3584} & 0.0997\%\\
          runtime & 12h& 20.4s & 0.0472\%\\
          memory & 41.4TiB & 1.69GiB & 0.0039\%\\
          \hline \\
         \textit{Fokker-Planck}  & full rank (CPU) & DLRA (GPU) & \% of full rank \\
            \hline 
          complexity & \num{11200000} & \num{11348642.8400} &101.3272\%\\
          system size & \num{11200000} & \num{567716.0000} & 5.0689\%\\
          runtime & 2684s &32.9s & 1.2258\%\\
          memory & 2.84TiB & 2.71GiB & 0.0932\%\\
          \hline
        \end{tabular}
    \end{table}
   Note, that a fair comparison to the costs of the TOPAS MC reference computation is not as straightforward \cite{borgers_complexity_1998}. For $10^8$ histories (primary particles) the computation time on the same workstation using up to 48 nodes is still in the ballpark of several hours. However, due to the different numerical discretization methods and included physics, the theoretical complexities are difficult to compare at a constant error. Further, the actual runtime can be almost arbitrarily reduced with more computational power due to the embarrassingly parallel nature of Monte Carlo methods. Thus, we give these values merely as a frame of reference.

\section{Discussion}\label{sec:discussion}
Our numerical results have demonstrated that the dynamical low-rank approximation can reduce the computational costs and memory associated with three-dimensional dose calculations in proton therapy by several orders of magnitude, both in theory and practice. In particular, we see that protons exhibit very low ranks even compared to other charged particles \cite{kusch_robust_2023}. This can likely be attributed to the structure of their energy deposition curve, i.e., the dynamics are mostly dominated by advection, followed by close to diffusive behavior for low energies. While ranks expectedly increase in heterogeneous materials, the fact that they do not increase significantly for several beam sources is promising for a use in treatment planning. 

We investigated two variations of the transport model: the Boltzmann equation and a Fokker-Planck approximation.  Overall, we observe slightly better results using the Boltzmann solver, both in terms of agreement with Monte Carlo, stability and ranks. Our comparative study between the Boltzmann and Fokker-Planck models reveals that, despite the Fokker-Planck model's simpler collision operator, the Boltzmann model provides more accurate results with lower computational ranks in our high-resolution framework. This suggests that the additional physical fidelity of the Boltzmann operator is crucial and is not a bottleneck when combined with DLRA. For both approaches, the low-rank solution shows near perfect agreement with the equivalent full-rank solution. Comparisons to the Monte Carlo reference show a need for high spatial and angular resolution, as solution quality deteriorates with coarser discretization. Yet, being able to relatively cheaply compute the dose distribution on such fine grids speaks for the efficacy of the dynamical low-rank approximation. The choice of spatial and angular discretizations however needs to be reevaluated, also with view to ``unphysical" behavior, i.e., negative regions in the solution. Further, the Fokker-Planck method heavily relies on the implicit scattering update and the correction described in \cite{landesman_angular_1989,morel_fokker-planck_1981} to avoid instabilities and artifacts at the material boundary, respectively. The correction however tends to overestimate the scattering moments \cite{drumm_analysis_2007} which could be a cause of the differences to Monte Carlo that we are seeing. Thus, future work could investigate whether the use of other, higher-order or locally refined spatial and angular discretizations can reduce the required expansion order and eliminate oscillations. Further, since protons undergo more physical interactions at low energies and thus energy loss is non-linear, a proportional, non-uniform energy/pseudo-time discretization could also be of interest to better capture rapid flux changes.

The currently used model for particle interactions still makes several simplifying assumptions. While the use of the continuous slowing down approximation and material composition model is common also in state-of-the-art Monte Carlo codes \cite{huang_validation_2018}, incorporating nuclear interactions as well as absorption is necessary in order to achieve clinical precision \cite{newhauser_physics_2015}. Note, that we do not see the full impact of neglecting these effects in this work as we compare to a Monte Carlo method which is also limited to electromagnetic interactions. Fortunately, adding nuclear interactions or absorption does not majorly affect the numerical method described here and is more a matter of accurate physical modeling or data (which can be challenging in itself). Based on the good agreement we achieve with the Monte Carlo reference and similar observations made in literature \cite{lathouwers_deterministic_2023}, neglecting straggling in the collided part of the equation appears to be a reasonable approximation. However, future work could also investigate strategies for handling a straggling term within the dynamical low-rank approximation of the collided equation.

\section{Conclusion}\label{sec:conclusion}
In this work, we demonstrated that the computational costs and memory associated with deterministic solvers for a proton transport problem can be significantly reduced using the dynamical low-rank method. The chosen ranks are low, making a much higher spatial and angular resolution computationally feasible at little additional costs. We show very good agreement with a full-rank solver based on the same physical models and discretizations. Further, our results agree reasonably well with a TOPAS MC reference. Future work could explore the benefits of higher-order or refined discretizations, include more accurate physics models and extend the framework to full radiation therapy treatment plans.
\section*{ACKNOWLEDGEMENTS}
Pia Stammer received funding from the German National Academy of Sciences Leopoldina for the project underlying this article, under grant number LPDS 2024-03.
\section*{CRediT author statement}
\textbf{Pia Stammer:} Conceptualization, Methodology, Software, Validation, Formal analysis, Investigation, Writing - Original Draft, Writing - Review \& Editing, Visualization \\
\textbf{Niklas Wahl:} Software, Writing - Review \& Editing \\
\textbf{Jonas Kusch:} Conceptualization, Methodology, Software, Writing - Original Draft, Writing - Review \& Editing \\
\textbf{Danny Lathouwers:} Methodology, Resources, Writing - Original Draft, Writing - Review \& Editing
\bibliographystyle{abbrv}
\bibliography{references}
\end{document}